\pdfoutput=1
\let\ORIlabel\label
\let\ORIrefstepcounter\refstepcounter
\AddToHook{package/hyperref/before}{\let\label\ORIlabel
   \let\refstepcounter\ORIrefstepcounter
}

\documentclass[onefignum,onetabnum]{siamart171218}

\usepackage{amsmath,amssymb,amsfonts,amsopn}
\usepackage{mathtools}
\usepackage{bm}
\usepackage{braket,amsopn}
\usepackage{array}

\usepackage{graphicx}

\ifpdf
  \DeclareGraphicsExtensions{.pdf,.png,.jpg,.eps}
\else
  \DeclareGraphicsExtensions{.eps}
\fi

\usepackage{epstopdf}

\makeatletter
\AtBeginDocument{\setkeys{Gin}{draft=false}\Gin@draftfalse
}
\makeatother

\usepackage{tikz}
\usetikzlibrary{positioning}

\usepackage{extarrows}
\usepackage{enumitem}
\usepackage{needspace}
\usepackage{placeins}

\usepackage{silence}
\usepackage{subcaption}
\usepackage{microtype}

\AtBeginDocument{
  \setlength{\abovedisplayskip}{6pt plus 2pt minus 2pt}
  \setlength{\belowdisplayskip}{6pt plus 2pt minus 2pt}
  \setlength{\abovedisplayshortskip}{4pt plus 2pt minus 2pt}
  \setlength{\belowdisplayshortskip}{4pt plus 2pt minus 2pt}
}

\newcommand{\real}{\mathbb{R}}

\newcommand{\fro}{\mathrm{F}}

\newcommand{\rank}{\mathrm{rank}}
\newcommand{\tr}{\mathrm{tr}}

\newcommand\vm{\mathbf{m}}

\newcommand\vx{\mathbf{x}}

\newcommand\ma{\mathbf{A}}
\newcommand\mb{\mathbf{B}}

\newcommand\md{\mathbf{D}}
\newcommand\me{\mathbf{E}}
\newcommand\mf{\mathbf{F}}
\newcommand\mg{\mathbf{G}}
\newcommand\mh{\mathbf{H}}
\newcommand\mi{\mathbf{I}}

\newcommand\mk{\mathbf{K}}
\newcommand\ml{\mathbf{L}}
\newcommand\mm{\mathbf{M}}

\newcommand\mq{\mathbf{Q}}
\newcommand\mP{\mathbf{P}}
\newcommand\mr{\mathbf{R}}
\newcommand\ms{\mathbf{S}}

\newcommand\mU{\mathbf{U}}
\newcommand\mv{\mathbf{V}}

\newcommand\mx{\mathbf{X}}
\newcommand\my{\mathbf{Y}}

\newcommand\mOmega{\mathbf{\Omega}}
\newcommand\mSigma{\mathbf{\Sigma}} 
\newcommand\mDelta{\mathbf{\Delta}}

\newcommand\mPsi{\mathbf{\Psi}}
\newcommand\mPi{\mathbf{\Pi}} 

\usepackage{cite}
\usepackage{booktabs}
\usepackage{url}
\usepackage{xcolor}
\usepackage{etoolbox} 

\usepackage{siunitx}
\usepackage{multirow}
\usepackage{colortbl}
\newcommand{\shadecell}{\cellcolor{black!10}} 



\usepackage{pgfplots}

\usepackage{algorithmic}
\makeatother

\Crefname{ALC@unique}{Line}{Lines}

\newsiamremark{remark}{Remark}
\crefname{remark}{Remark}{Remarks}

\headers{Adaptive Matrix-Free HODLR Approximation}{H. Bhatt, A. Alexanderian, C. Chen}

\title{Adaptive Matrix-Free HODLR Approximation with Global Error Control}

\author{Harshit Bhatt\thanks{Department of Mathematics, North Carolina State University 
  (\email{hbhatt@ncsu.edu}, \email{aalexan3@ncsu.edu}, \email{cchen49@ncsu.edu}).}
\and Alen Alexanderian\footnotemark[1]
\and Chao Chen\footnotemark[1]}

\makeatletter
\newcommand*{\addFileDependency}[1]{\typeout{(#1)}\@addtofilelist{#1}\IfFileExists{#1}{}{\typeout{No file #1.}}}
\makeatother

\pgfplotsset{compat=1.18}

\begin{document}
\maketitle

\begin{abstract}
Many dense matrices in scientific computing, such as Schur complements and Hessians, are too expensive to form but can be applied to vectors efficiently. In many cases, the off-diagonal blocks of these matrices have low numerical rank, which a hierarchically off-diagonal low-rank (HODLR) approximation exploits. Existing methods that build this approximation from matrix-vector products require the ranks in advance or do not bound the error of the full approximation. 
We present an adaptive randomized algorithm that requires only a relative tolerance. Our method compresses the off-diagonal blocks from the largest to the smallest and enlarges each approximation until an error test on new random vectors meets the tolerance. 
When the ranks are bounded, the number of matrix-vector products grows logarithmically with the matrix size. Because the largest blocks are compressed first, their errors enter the random samples of the smaller blocks. We prove that these errors do not accumulate: when the tests use enough random vectors, the relative Frobenius-norm error of the full approximation is, with high probability, at most a constant times the tolerance, and the constant is independent of the matrix size. Numerical results on Schur complements and an inverse-problem Hessian demonstrate the robustness, efficiency, and asymptotic scaling of our algorithm and provide strong evidence for our theory: for tolerances from $10^{-1}$ to $10^{-14}$, the errors are at most $1.1$ times the tolerance. \end{abstract}



\section{Introduction}
High-dimensional dense matrices arise throughout scientific computing and machine learning, including discretized integral operators, kernel and covariance matrices, Schur complements, and Hessians. Storing a general $n\times n$ matrix requires $\mathcal{O}(n^2)$ memory, and standard dense factorizations require $\mathcal{O}(n^3)$ work. Hierarchical matrix methods reduce these costs by recursively partitioning the matrix and approximating compressible off-diagonal blocks with low-rank factors. When the numerical ranks remain bounded, these representations support storage, matrix-vector products (matvecs), factorizations, and linear solves at nearly linear cost~\cite{hackbuschHierarchicalMatricesAlgorithms2015}.

Constructing a hierarchical representation is challenging when individual entries are unavailable or expensive to query. In the \emph{matrix-free} setting considered here, we can apply the matrix and its transpose to vectors without accessing individual matrix entries. Our objective is to construct a hierarchically off-diagonal low-rank (HODLR) approximation with controlled relative Frobenius-norm error. Existing matrix-free HODLR constructions with an error analysis require ranks or rank bounds in advance. The prescribed ranks or bounds determine the number of random vectors used to approximate the off-diagonal blocks. Our method chooses this number and the ranks itself: it tests the error of the current approximation on new random vectors and enlarges the approximation until the error test is satisfied. In the terms used for low-rank approximation~\cite{halkoFindingStructureRandomness2011,yuEfficientRandomizedAlgorithms2018}, these constructions address the fixed-rank problem, where the ranks are given. We address the fixed-precision problem, in which the accuracy is specified.

HODLR uses a simple recursive partition in which every off-diagonal block is compressed (weak admissibility), and it does not require the nested bases of the hierarchically semiseparable (HSS) and $\mathcal{H}^2$ formats. It is therefore useful when stronger geometric assumptions or basis-sharing constraints are unavailable. Existing software, including HODLRlib~\cite{ambikasaranHODLRlib2019}, the MATLAB hm-toolbox~\cite{masseiHMToolbox2020}, and a GPU-accelerated factorization and solver~\cite{chenSolvingLinearSystems2022}, provides arithmetic routines for the resulting representations.

Two applications motivate our construction. Eliminating variables in sparse direct solvers produces dense Schur complements and frontal matrices~\cite{davisSurveyDirectMethods2016}. Schur-complement actions can often be computed using sparse operations and subdomain solves without forming the dense matrix. Similarly, in PDE-constrained inverse problems, a data-misfit Hessian-vector product requires forward and adjoint PDE solves, while individual entries are not readily accessible~\cite{ambartsumyanHierarchicalMatrixApproximations2020,bui-thanhComputationalFrameworkInfiniteDimensional2013}. Hierarchical approximations exploit localized off-diagonal structure when a global low-rank approximation is insufficient~\cite{ambartsumyanHierarchicalMatrixApproximations2020,hartlandHierarchicalOffdiagonalLowrank2023}. Both applications call for accuracy control with economical operator access.

\subsection{Related work}
\paragraph{Entry-based construction}
Many methods construct an HODLR approximation from selected matrix entries and choose the ranks adaptively, for example, adaptive cross approximation and boundary-distance low-rank approximation~\cite{aminfarFastBlockLowRank2016,ambikasaranHODLRlib2019,masseiHMToolbox2020}. They cannot be applied in the matrix-free setting, where individual entries are not available.

\paragraph{Matrix-free construction and error analysis}
Lin et al.~\cite{linFastConstructionHierarchical2011} introduced top-down peeling for hierarchical matrices. Peeling compresses the blocks one level of the recursive partition at a time, from the largest blocks to the smallest, and subtracts each computed approximation before it samples the next level. Martinsson~\cite{martinssonCompressingRankStructuredMatrices2016} developed a matrix-free randomized HODLR construction of this type. These methods assume prescribed off-diagonal rank bounds. Hartland et al.~\cite{hartlandHierarchicalOffdiagonalLowrank2023} report using a rank-adaptive variant for data-misfit Hessians, with a tolerance for each block, and note that their error estimate neglects the error incurred in peeling. Chen et al.~\cite{chenNearoptimalHierarchicalMatrix2025} use generalized Nystr\"om approximation and randomly perforated sketches to obtain a near-optimal expected squared Frobenius-norm error relative to the best rank-$k$ HODLR approximation. Their formulation prescribes $k$, whereas ours selects the ranks from a given tolerance. Adaptive construction from matvecs and selected entries has also been studied for HSS and $\mathcal{H}^2$ formats~\cite{gormanRobustAccurateStopping2019,yanivAdaptiveJLHSS2025,boukaramAdaptiveSketchingH22025}. These works do not analyze how the local stopping tests control the error of the full matrix.

The accumulation of errors across levels is a natural concern in peeling~\cite{linFastConstructionHierarchical2011}. Chen et al.~\cite{chenNearoptimalHierarchicalMatrix2025} show that rank-truncated randomized-SVD peeling can amplify errors exponentially in the depth, and control this accumulation in their modified fixed-rank algorithm. Carson et al.~\cite{carsonMixedPrecisionHODLR2025} analyze HODLR approximation and mixed-precision storage when the blocks are available for independent compression. Our analysis instead accounts for the errors of earlier levels that enter the samples of later levels in adaptive matrix-free peeling.

\subsection{Contributions}
We show that an HODLR approximation can be computed using matvecs to a prescribed accuracy without prior knowledge of the ranks. Given a relative tolerance $\varepsilon$, the algorithm determines the ranks adaptively, and the relative Frobenius-norm error of the full matrix is, with high probability, bounded by a constant multiple of $\varepsilon$.

\begin{itemize}
    \item \textbf{The tolerance determines the ranks (\S\ref{sec:adaptivepeeling}, \S\ref{sec:recompression}).} Our algorithm, adaptive peeling, is an adaptive version of the peeling method. It compresses the off-diagonal blocks from the largest to the smallest sizes. It enlarges the approximation of all blocks of one size by $b$ columns at a time and tests the error on new random vectors until a randomized estimate of the relative error is below $\varepsilon$. No rank bound is needed, and blocks of different sizes receive different ranks. We strongly recommend a recompression step after construction, which removes the columns that the tolerance does not require. The block size $b$ can therefore be chosen for reliable tests and efficient matrix operations, and not for the ranks. When the ranks are bounded by a constant $r$, the construction uses $\mathcal{O}(r\log n)$ matvecs for an $n\times n$ matrix, and the approximation requires $\mathcal{O}(rn\log n)$ storage, as in peeling with prescribed ranks. The dense diagonal blocks are recovered with Rademacher variables in place of the stacked identity matrices of fixed-rank peeling~\cite[\S4]{martinssonCompressingRankStructuredMatrices2016}, which keeps the recovery error at the level of the off-diagonal error (\Cref{lem:diagonal_error}).

    \item \textbf{The tolerance bounds the error of the whole matrix (\S\ref{sec:global_error_analysis}).} With high probability, the returned approximation $\widetilde{\mh}$ satisfies $\|\mh-\widetilde{\mh}\|_\fro\leq C\varepsilon\|\mh\|_\fro$, where $C$ does not depend on the matrix size (\Cref{thm:global_error}). The difficulty is that the samples for smaller blocks contain the errors already made on larger blocks. We show that each test measures these inherited errors together with the new compression error. A test that passes therefore accounts for both, and the errors do not accumulate. Inherited errors can raise the ranks or leave larger dense diagonal blocks, which raises the cost, but they cannot raise the error above the bound. The analysis also gives a number of random vectors per test that is sufficient.

    \item \textbf{The experiments demonstrate the algorithm and support the analysis (\S\ref{sec:results}).} We apply adaptive peeling to representative problems: Schur complements that arise in sparse direct solvers and the data-misfit Hessian of a PDE-constrained inverse problem. The results demonstrate the scaling, efficiency, and robustness of the algorithm. They also provide strong numerical evidence for the analysis, including the role of the block size, the success probability of the sampled tests, and the effect of the random samples.
\end{itemize}
\vspace{.5em}

\paragraph{Notation}
Matrices are denoted by bold uppercase and vectors by bold lowercase letters. Given index sets $I,J\subseteq\{1,\ldots,n\}$, we write $\ma(I,J)$ for the submatrix of $\ma$ with rows indexed by $I$ and columns indexed by $J$, and a colon denotes the full index range, so that $\ma(I,:)$ selects the rows in $I$ and $\ma(:,J)$ the columns in $J$. Throughout, $\|\cdot\|_\fro$ denotes the Frobenius norm, and $\|\cdot\|_2$ denotes the spectral norm for matrices and the Euclidean norm for vectors. The notation $\operatorname{diag}(\vx)$ forms a diagonal matrix from a vector.

\section{Background}
\label{sec:background}
In this section, we review the HODLR matrix structure and Gaussian estimates of the matrix Frobenius norm.

\subsection{HODLR matrices}
\label{sec:hodlr}

To define an HODLR matrix, we first describe the binary-tree partitioning that is used to identify subblocks of the matrix. An HODLR representation with $L$ levels then stores the diagonal blocks of the finest level densely and approximates the off-diagonal blocks that couple the two children of every tree node by low-rank factorizations~\cite{ambikasaranDarveHODLR2013,martinssonCompressingRankStructuredMatrices2016}.

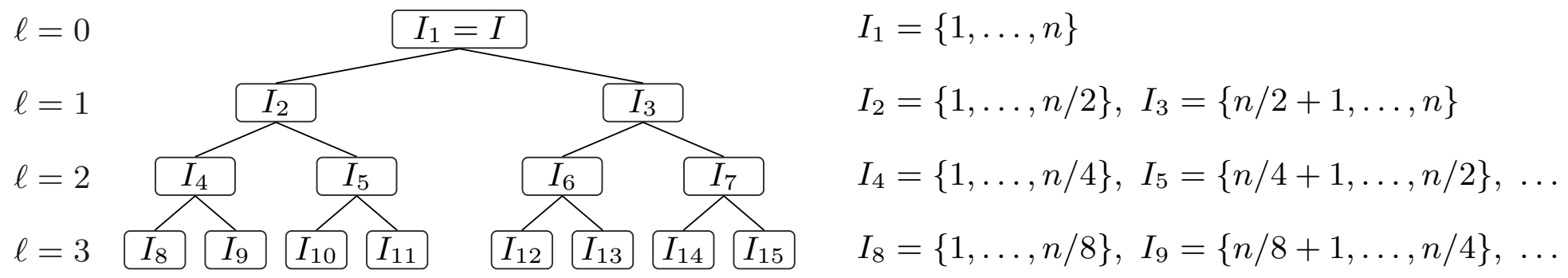
\begin{figure}
\centering
\resizebox{1\linewidth}{!}{\begin{tikzpicture}[
    x=0.35cm,
    y=0.52cm,
    tnode/.style={
        draw=black!85,
rounded corners=1.5pt,
        minimum width=0.76cm,
        minimum height=0.36cm,
        inner sep=1.2pt,
        font= \small
        },
    root/.style={
        tnode,
minimum width=1.28cm
    },
    leaf/.style={
        tnode,
minimum width=0.58cm,
        inner sep=1.0pt
    },
    branch/.style={
line width=0.4pt
    },
    lev/.style={
        font= \small,
        anchor=east,
        text=black!90
    },
    info/.style={
rounded corners=2pt,
        inner xsep=3pt,
        inner ysep=3pt,
        font=\small,
        anchor=west,
        align=left,
        text width=6.7cm
    }
]

\node[root] (I1) at (0,0) {$I_1=I$};

\node[tnode] (I2) at (-5.0,-1.35) {$I_2$};
\node[tnode] (I3) at ( 5.0,-1.35) {$I_3$};

\node[tnode] (I4) at (-7.2,-2.70) {$I_4$};
\node[tnode] (I5) at (-2.8,-2.70) {$I_5$};
\node[tnode] (I6) at ( 2.8,-2.70) {$I_6$};
\node[tnode] (I7) at ( 7.2,-2.70) {$I_7$};

\node[leaf] (I8)  at (-8.3,-4.05) {$I_8$};
\node[leaf] (I9)  at (-6.1,-4.05) {$I_9$};
\node[leaf] (I10) at (-3.9,-4.05) {$I_{10}$};
\node[leaf] (I11) at (-1.7,-4.05) {$I_{11}$};
\node[leaf] (I12) at ( 1.7,-4.05) {$I_{12}$};
\node[leaf] (I13) at ( 3.9,-4.05) {$I_{13}$};
\node[leaf] (I14) at ( 6.1,-4.05) {$I_{14}$};
\node[leaf] (I15) at ( 8.3,-4.05) {$I_{15}$};

\draw[branch] (I1.south) -- (I2.north);
\draw[branch] (I1.south) -- (I3.north);

\draw[branch] (I2.south) -- (I4.north);
\draw[branch] (I2.south) -- (I5.north);
\draw[branch] (I3.south) -- (I6.north);
\draw[branch] (I3.south) -- (I7.north);

\draw[branch] (I4.south) -- (I8.north);
\draw[branch] (I4.south) -- (I9.north);
\draw[branch] (I5.south) -- (I10.north);
\draw[branch] (I5.south) -- (I11.north);
\draw[branch] (I6.south) -- (I12.north);
\draw[branch] (I6.south) -- (I13.north);
\draw[branch] (I7.south) -- (I14.north);
\draw[branch] (I7.south) -- (I15.north);

\node[lev] at (-9.7,  0.00) {$\ell=0$};
\node[lev] at (-9.7, -1.35) {$\ell=1$};
\node[lev] at (-9.7, -2.70) {$\ell=2$};
\node[lev] at (-9.7, -4.05) {$\ell=3$};

\node[info] at (10.5, 0.00) {
    $\displaystyle I_1=\{1,\ldots,n\}$
};

\node[info] at (10.5,-1.35) {
    $\displaystyle
    I_2=\{1,\ldots,n/2\},\
    I_3=\{n/2+1,\ldots,n\}$
};

\node[info] at (10.5,-2.70) {
    $\displaystyle
    I_4=\{1,\ldots,n/4\},\
    I_5=\{n/4+1,\ldots,n/2\},\ \ldots$
    };

\node[info] at (10.5,-4.05) {
    $\displaystyle
    I_8=\{1,\ldots,n/8\},\
    I_9=\{n/8+1,\ldots,n/4\},\ \ldots$
};

\end{tikzpicture}}
\caption{A complete balanced binary partition tree of depth $L=3$, with the
root at level $0$ and the leaves at level $3$.}
\label{fig:binary_hodlr_tree}
\end{figure}

\paragraph{Binary-tree partitioning}
Let $I=\{1,\ldots,n\}$, and apply the same partition to the row and column indices of $\mh\in\real^{n\times n}$. Throughout the paper, we use a complete balanced tree of depth $L$ and assume for simplicity that $2^L$ divides $n$. The root node has index set $I_1=I$ at level $0$. We number the nodes breadth first: the left and right children of node $p$ are $2p$ and $2p+1$, respectively, with
\[
    I_p = I_{2p} \,\cup\, I_{2p+1},
    \qquad I_{2p} \cap I_{2p+1} = \emptyset.
\]
The children are \emph{siblings}, and nodes at level $L$ are \emph{leaves}. Level $\ell$ contains nodes $2^\ell,\ldots,2^{\ell+1}-1$, each with $n/2^\ell$ indices. We further require that every $I_\tau$ consists of \emph{consecutive} indices, so that a node at level $\ell$ has $I_\tau=\{p_\tau+1,\ldots,p_\tau+n/2^{\ell}\}$ for some offset $p_\tau$, and its two children take the two halves of this range; see \Cref{fig:binary_hodlr_tree}. Contiguity ensures that sibling interactions are contiguous submatrices of $\mh$; these are the off-diagonal blocks compressed by the HODLR format. In practice, this requirement is met by applying a suitable permutation to $\mh$ before the index partitioning, so that indices that are close in the underlying geometry are numbered consecutively. 

For operators with a geometric origin, we use a kd-tree ordering. Starting from all points, each group is split at the median of the coordinate with the largest span, recursively, until every group has at most a prescribed number of points. \Cref{fig:kdtree_levels} shows three levels of this bisection with the node numbers of \Cref{fig:binary_hodlr_tree}; the 3D Schur complement and the Hessian of \S\ref{sec:results} are ordered this way. 
 The algorithm and the analysis do not rely on equal splits and extend to trees with nonuniform splits; we use the balanced tree to keep the notation simple.

\begin{figure}[ht]
\centering
\newcommand{\kdpanel}[2]{\begin{tikzpicture}[x=1cm,y=1cm,baseline=(current bounding box.south)]
  \foreach \i in {0,...,9} {
    \foreach \j in {0,...,9} {
      \fill[gray!70] (0.21*\i+0.105,0.21*\j+0.105) circle (0.4pt);
    }
  }
  \draw[thin] (0,0) rectangle (2.1,2.1);
  #1
  \node[font=\footnotesize] at (1.05,-0.3) {#2};
\end{tikzpicture}}
\newcommand{\kdlabel}[3]{\node[font=\footnotesize,fill=white,inner sep=1pt] at (#1,#2) {#3};}
\kdpanel{\draw[line width=0.9pt] (1.05,0) -- (1.05,2.1);
  \kdlabel{0.525}{1.05}{2}\kdlabel{1.575}{1.05}{3}}{level $1$: nodes $2$, $3$}\hspace{3em}
\kdpanel{\draw[line width=0.9pt] (1.05,0) -- (1.05,2.1);
  \draw[line width=0.6pt] (0,1.05) -- (2.1,1.05);
  \kdlabel{0.525}{0.525}{4}\kdlabel{0.525}{1.575}{5}\kdlabel{1.575}{0.525}{6}\kdlabel{1.575}{1.575}{7}}{level $2$: nodes $4$--$7$}\hspace{3em}
\kdpanel{\draw[line width=0.9pt] (1.05,0) -- (1.05,2.1);
  \draw[line width=0.6pt] (0,1.05) -- (2.1,1.05);
  \draw[line width=0.35pt] (0.525,0) -- (0.525,2.1);
  \draw[line width=0.35pt] (1.575,0) -- (1.575,2.1);
  \kdlabel{0.2625}{0.525}{8}\kdlabel{0.7875}{0.525}{9}\kdlabel{0.2625}{1.575}{10}\kdlabel{0.7875}{1.575}{11}\kdlabel{1.3125}{0.525}{12}\kdlabel{1.8375}{0.525}{13}\kdlabel{1.3125}{1.575}{14}\kdlabel{1.8375}{1.575}{15}}{level $3$: leaves $8$--$15$}
 \caption{The kd-tree ordering on a square point set. Each node of \Cref{fig:binary_hodlr_tree} is a spatial cluster, split at the median of its longer side into its two children. Left: nodes $2$ and $3$ at level $1$. Middle: nodes $4$--$7$ at level $2$. Right: the leaves $8$--$15$ at level $3$. Numbering the points cluster by cluster in this order gives the consecutive index sets $I_\tau$.}
\label{fig:kdtree_levels}
\end{figure}
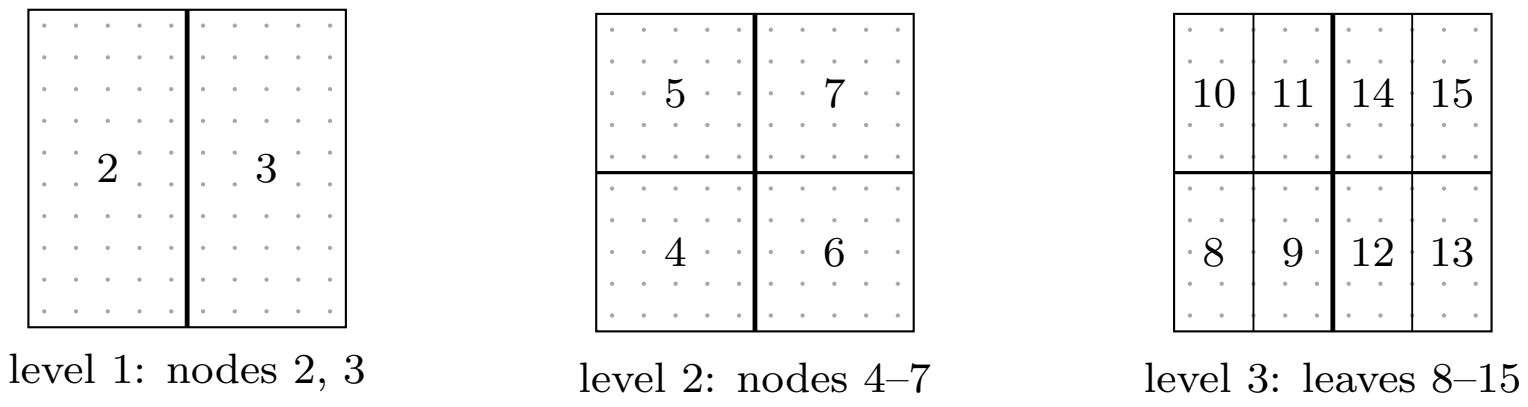

\begin{figure}
    \centering
{\definecolor{blockorange}{rgb}{0.95,0.64,0.17}
\colorlet{blockedge}{blockorange!75!black}

\def\matrixsize{1.9cm}
\def\matrixgap{0.4}
\def\signcolor{black}

\pgfmathsetmacro{\step}{1+\matrixgap}
\pgfmathsetmacro{\signpos}{1+\matrixgap/2}

\begin{tikzpicture}[x=\matrixsize,y=\matrixsize,
block/.style={fill=blockorange,draw=blockedge,line width=.2pt},
    matrix frame/.style={draw=black,line width=.4pt},
    block label/.style={font=\scriptsize,inner sep=0pt}]

\fill[blockorange] (0,0) rectangle (1,1);
\draw[matrix frame] (0,0) rectangle (1,1);
\node at (.5,1.15) {$\mh$};

\node[\signcolor] at (\signpos,.5) {$=$};
\foreach \j in {1,2,3}{
    \node[\signcolor] at ({\j*\step+\signpos},.5) {$+$};
}

\begin{scope}[shift={(\step,0)}]
    \foreach \j in {0,...,7}{
        \path[block] ({\j/8},{1-(\j+1)/8})
            rectangle ({(\j+1)/8},{1-\j/8});
    }
    \draw[matrix frame] (0,0) rectangle (1,1);
    \node at (.5,1.15) {$\md_3$};
\end{scope}

\begin{scope}[shift={({2*\step},0)}]
    \path[block] (.5,.5) rectangle (1,1);
    \path[block] (0,0) rectangle (.5,.5);
    \draw[matrix frame] (0,0) rectangle (1,1);
    \node at (.5,1.15) {$\mh^{(1)}$};
    \node[block label] at (.75,.75) {$\mh_{23}$};
    \node[block label] at (.25,.25) {$\mh_{32}$};
\end{scope}

\begin{scope}[shift={({3*\step},0)}]
    \path[block] (.25,.75) rectangle (.5,1);
    \path[block] (0,.5) rectangle (.25,.75);
    \path[block] (.75,.25) rectangle (1,.5);
    \path[block] (.5,0) rectangle (.75,.25);
    \draw[matrix frame] (0,0) rectangle (1,1);
    \node at (.5,1.15) {$\mh^{(2)}$};
    \node[block label] at (.375,.875) {$\mh_{45}$};
    \node[block label] at (.125,.625) {$\mh_{54}$};
    \node[block label] at (.875,.375) {$\mh_{67}$};
    \node[block label] at (.625,.125) {$\mh_{76}$};
\end{scope}

\begin{scope}[shift={({4*\step},0)}]
    \foreach \j in {0,...,3}{
        \pgfmathsetmacro{\offset}{\j/4}
        \path[block] ({\offset+.125},{1-\offset-.125})
            rectangle ({\offset+.25},{1-\offset});
        \path[block] (\offset,{1-\offset-.25})
            rectangle ({\offset+.125},{1-\offset-.125});
    }
    \draw[matrix frame] (0,0) rectangle (1,1);
    \node at (.5,1.15) {$\mh^{(3)}$};
\end{scope}

\end{tikzpicture}

 }
    \caption{Exact block decomposition induced by the binary partition tree of depth $L=3$. The labels identify sibling blocks $\mh_{ab}=\mh(I_a,I_b)$ using the node indices in \Cref{fig:binary_hodlr_tree}; reversing the indices selects the opposite directed block.}
 \label{fig:exact_block_decomposition}
\end{figure}
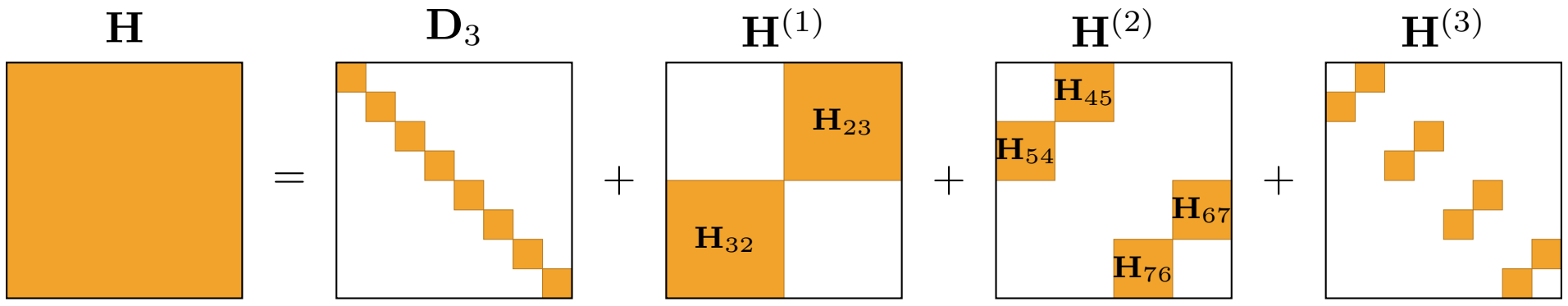

Write $\mh_{ab}:=\mh(I_a,I_b)$. At level $\ell$, the sibling blocks are $\mh_{2p,2p+1}$ and $\mh_{2p+1,2p}$ for parents $p=2^{\ell-1},\ldots,2^\ell-1$. Let $\mh^{(\ell)}\in\real^{n\times n}$ equal $\mh$ on these blocks and zero elsewhere, and let $\md_L$ equal $\mh$ on the leaf diagonal blocks and zero elsewhere. The tree then induces the exact decomposition
\begin{equation}
\label{eq:hodlr_exact_decomposition}
    \mh=\md_L+\sum_{\ell=1}^{L}\mh^{(\ell)}.
\end{equation}
For example, in \Cref{fig:exact_block_decomposition}, $\mh^{(1)}$ contains the interaction between $I_2$ and $I_3$, whereas $\mh^{(3)}$ contains the finest sibling interactions. The matrices $\md_L,\mh^{(1)},\ldots,\mh^{(L)}$ have pairwise disjoint supports. Consequently,
\begin{equation}
\label{eq:hodlr_energy_decomposition}
    \|\mh\|_\fro^2
    =\|\md_L\|_\fro^2+\sum_{\ell=1}^{L}\|\mh^{(\ell)}\|_\fro^2.
\end{equation}

The decomposition \eqref{eq:hodlr_exact_decomposition} is exact for every matrix; it becomes useful only when the sibling off-diagonal blocks are compressible.

\begin{definition}[HODLR matrix]
\label{def:hodlr}
Let $k_1,\ldots,k_L$ be nonnegative integers. A matrix $\mh\in\real^{n\times n}$ is an \emph{HODLR matrix of ranks $(k_1,\ldots,k_L)$} with respect to a binary tree partitioning if, for every level $\ell=1,\ldots,L$ and parent $p=2^{\ell-1},\ldots,2^\ell-1$,
\[
    \rank(\mh_{2p,2p+1})\le k_\ell,
    \qquad \rank(\mh_{2p+1,2p})\le k_\ell.
\]
For a generic sibling pair, write its node indices as $\alpha$ and $\beta$. The corresponding \emph{HODLR representation} stores factors $\mU_{\alpha\beta},\mv_{\alpha\beta}\in\real^{(n/2^\ell)\times k_\ell}$ satisfying $\mh_{\alpha\beta}=\mU_{\alpha\beta}\mv_{\alpha\beta}^\top$ for every sibling pair, together with the leaf diagonal blocks $\md_L$ in dense form.
\end{definition}

Matrices arising in applications are seldom HODLR in this exact sense. Instead, the singular values of their sibling off-diagonal blocks decay rapidly, so that each block is well approximated by a low-rank factorization with a rank much smaller than the size of the block; the Schur complements and the data-misfit Hessian of \S\ref{sec:results} are of this type. These numerical ranks are not known a priori and may vary across blocks and levels. Replacing each sibling block in \eqref{eq:hodlr_exact_decomposition} by a low-rank approximation gives an HODLR approximation 
\begin{equation}
\label{eq:HODLR_approximation}
\widetilde{\mh}=\widetilde{\md}_L+\sum_{\ell=1}^{L}\widetilde{\mh}^{(\ell)},
\end{equation}
where the leaf diagonal blocks in $\widetilde{\md}_L$ are stored densely. If $\mh$ is symmetric, then $\mh_{\beta\alpha}=\mh_{\alpha\beta}^\top$, so only one factorization per unordered sibling pair is needed. A nonsymmetric matrix requires both directed blocks.

\subsection{Gaussian error estimates}
\label{sec:gaussian_error_estimates}
\label{sec:hutchinson}
To choose the ranks, the algorithm must know the error of each approximation $\widetilde{\mh}^{(\ell)}$.  Since the blocks of $\mh$ are not explicitly available, this error cannot be computed directly. We estimate it from matvecs. Given an approximation $\widetilde{\ma}$ to $\ma$, we estimate its relative Frobenius error using a Gaussian test matrix $\mOmega$:
\begin{equation}
\label{eq:gaussian_error_indicator}
\widehat e
:=\frac{\|(\ma-\widetilde{\ma})\mOmega\|_\fro}{\|\ma\mOmega\|_\fro}.
\end{equation}

This ratio can be evaluated using matvecs. Throughout the paper, a hat marks a quantity computed from random samples. The following corollary uses the Gaussian norm estimate of Gratton and Titley-Peloquin~\cite[Theorem~3.1]{grattonImprovedBoundsSmallSample2018} for the numerator and denominator. The required number of random vectors depends on stable ranks. For nonzero $\mx$, the \emph{stable rank} is $\|\mx\|_\fro^2/\|\mx\|_2^2$ and lies between $1$ and the rank of $\mx$.

\begin{corollary}[Gaussian relative-error estimate]
\label{thm:stability_hutchinson}
Let $\ma,\widetilde{\ma}\in\real^{m\times n}$ be fixed, with $\|\ma\|_\fro>0$, and define $\mDelta:=\ma-\widetilde{\ma}\ne0$. Let $\mOmega\in\real^{n\times b}$, $b\leq n$, have independent standard normal entries. Let $\rho_*$ be the smaller of the stable ranks of $\ma$ and $\mDelta$,
\[
\rho_*:=\min\left\{
\frac{\|\ma\|_\fro^2}{\|\ma\|_2^2},
\frac{\|\mDelta\|_\fro^2}{\|\mDelta\|_2^2}
\right\}.
\]
For any $\eta\in(0,1/2]$ and $\delta\in(0,1)$, set $c_\eta:=(1+\eta)/(1-\eta)$. If $b\ge 2\rho_*^{-1}\eta^{-2}\log(2/\delta)$, then, with probability at least $1-\delta$,
\[
\frac{\|\ma-\widetilde{\ma}\|_\fro}{\|\ma\|_\fro}
\le c_\eta\widehat e.
\]
\end{corollary}

\begin{proof}
Apply the lower tail bound of \cite[Theorem~3.1]{grattonImprovedBoundsSmallSample2018} to $\mDelta$ with factor $\tau=(1-\eta)^{-1}$ and the upper tail bound to $\ma$ with $\tau=1+\eta$. Since $(1-\tau^{-2})^2=\eta^2(2-\eta)^2\geq2\eta^2$ for the lower tail, each failure probability is at most $\exp(-b\rho_*\eta^2/2)\leq\delta/2$. Thus, with probability at least $1-\delta$,
\[
\|\mDelta\mOmega\|_\fro\geq\sqrt b(1-\eta)\|\mDelta\|_\fro,
\qquad
\|\ma\mOmega\|_\fro\leq\sqrt b(1+\eta)\|\ma\|_\fro.
\]
Taking the ratio gives the claimed bound.
\end{proof}

Consequently, if a test passes, that is, $\widehat e\le\varepsilon$, then $\|\ma-\widetilde{\ma}\|_\fro\le c_\eta\varepsilon\|\ma\|_\fro$ on the same event. In the rest of the paper, we use \Cref{thm:stability_hutchinson} with $\eta=1/2$. Then $c_\eta=3$, and the condition on the width is $b\geq8\rho_*^{-1}\log(2/\delta)$. A wider test fails less often and has a constant closer to $1$. In the algorithm of \S\ref{sec:methods}, however, $b$ is also the number of columns by which the ranks grow.

\section{Adaptive matrix-free HODLR approximation}
\label{sec:methods}

Peeling~\cite{linFastConstructionHierarchical2011,martinssonCompressingRankStructuredMatrices2016} constructs an HODLR approximation from matvecs when the ranks of the sibling blocks are prescribed. We develop an adaptive version of this method, which we call \emph{adaptive peeling}, for a matrix $\mh\in\real^{n\times n}$ available only through matvecs. Instead of ranks, the user prescribes a relative tolerance $\varepsilon\in(0,1)$, which is used as in \S\ref{sec:gaussian_error_estimates}. At every level $\ell$, the algorithm enlarges the approximation $\widetilde{\mh}^{(\ell)}$ of the sibling blocks $\mh^{(\ell)}$ until the error estimate $\widehat e$ of \eqref{eq:gaussian_error_indicator} satisfies $\widehat e\leq\varepsilon$. The construction then satisfies
\[
    \|\mh-\widetilde{\mh}\|_\fro
    \leq \mathcal{O}(\varepsilon)\|\mh\|_\fro
\]
with high probability under the sampling and recovery conditions of \Cref{thm:global_error}. \S\ref{sec:adaptivepeeling} describes how the algorithm turns the tolerance into ranks, \S\ref{sec:global_error_analysis} proves the bound, and \S\ref{sec:recompression} reduces the ranks after construction.

\subsection{Adaptive peeling}
\label{sec:adaptivepeeling}
We use the balanced partition of \S\ref{sec:hodlr}. At level $\ell$, let $\alpha=2p$ and $\beta=2p+1$ be the children of $p=2^{\ell-1},\ldots,2^\ell-1$. The level matrix $\widetilde{\mh}^{(\ell)}$ contains approximations $\widetilde{\mh}_{\alpha\beta}$ and $\widetilde{\mh}_{\beta\alpha}$ of these sibling blocks and is zero elsewhere.

Peeling compresses the sibling blocks one level at a time, from the coarsest level to the finest~\cite{linFastConstructionHierarchical2011,martinssonCompressingRankStructuredMatrices2016}. After a level is compressed, its approximation is subtracted, and the next level is compressed from matvecs with the remaining matrix. After $k$ compressed levels, this matrix is the residual
\begin{equation}
\label{eq:peeling_residual}
\mr^{(k)}:=\mh-\sum_{\ell=1}^{k}\widetilde{\mh}^{(\ell)},
\qquad k=0,\ldots,L,
\end{equation}
where $L$ is the number of compressed levels and $\mr^{(0)}=\mh$. The level-$(k+1)$ sibling blocks of $\mr^{(k)}$ are those of $\mh$, and its coarser sibling blocks hold the errors of the $k$ computed approximations. After the last level, the leaf diagonal blocks $\widetilde{\md}_L$ are recovered from $\mr^{(L)}$, and the returned approximation is $\widetilde{\mh}$ in \eqref{eq:HODLR_approximation}.

At each level, adaptive peeling enlarges the low-rank bases by $b$ columns at a time, as in the blocked randomized QB (randQB) algorithms~\cite{yuEfficientRandomizedAlgorithms2018}. After each enlargement, we draw new random vectors and compute the error estimate of \S\ref{sec:gaussian_error_estimates}, which requires only matvecs; this is the test that Smith et al.~\cite{smithAdaptiveMatrixFreeLowRank2026} use for a single matrix. If the estimate is at most $\varepsilon$, the level is \emph{accepted}: its approximation is subtracted, and the algorithm moves to the next level. Otherwise, the bases are enlarged again. The compressed levels are therefore the accepted levels, and we call their number $L$ the \emph{accepted depth}. It can be smaller than the intended depth $L_{\max}$ given to the algorithm, where $2^{L_{\max}}$ divides $n$, because the algorithm stops when the rank at a level reaches a cutoff; \Cref{rem:termination} states the rule.

For ease of presentation, we first explain the method for a symmetric matrix. The modification for a nonsymmetric matrix, which requires matvecs with both $\mh$ and $\mh^\top$, is given in \Cref{rem:nonsymmetric}.

\subsubsection{Level \texorpdfstring{$1$}{1} for a symmetric matrix: adaptive compression}
The first level contains two sibling blocks, $\mh_{23}$ and $\mh_{32}=\mh_{23}^\top$. We build a basis for the range of $\mh_{23}$ in increments of $b$ columns and index the increments by $i=1,2,\ldots$. Increment $i$ draws a fresh standard Gaussian matrix $\mg_{3;i}\in\real^{(n/2)\times b}$, where the subscript $3$ records the node whose columns are sampled, and $i$ the increment. To sample $\mh_{23}$ through a product with $\mh$, place $\mg_{3;i}$ on the rows in $I_3$ and zeros on the rows in $I_2$:
\[
\mOmega_i=
\begin{bmatrix}\mathbf 0\\\mg_{3;i}\end{bmatrix},
\qquad
\ms_i=\mh\mOmega_i=
\begin{bmatrix}\mh_{22}&\mh_{23}\\\mh_{32}&\mh_{33}\end{bmatrix}
\begin{bmatrix}\mathbf 0\\\mg_{3;i}\end{bmatrix}
=
\begin{bmatrix}\mh_{23}\mg_{3;i}\\\mh_{33}\mg_{3;i}\end{bmatrix}.
\]
The zero rows remove $\mh_{22}$ from the product, so the upper row block of $\ms_i$, $\my_{2;i}:=\mh_{23}\mg_{3;i}$, is a Gaussian sample of the range of $\mh_{23}$. The lower row block is not used at this level.

For $i=1$, a thin QR factorization~\cite[\S5.2]{golubMatrixComputations2013} of $\my_{2;1}$ gives a basis $\mq_2\in\real^{(n/2)\times b}$ with orthonormal columns. To express $\mh_{23}$ in this basis, stack $\mq_2$ with zeros and compute $\mb:=(\mh^\top\mq)^\top=\mq^\top\mh$. Specifically,
\[
\mq=\begin{bmatrix}\mq_2\\\mathbf 0\end{bmatrix},
\qquad
\mq^\top\mh=
\begin{bmatrix}\mq_2^\top&\mathbf 0\end{bmatrix}
\begin{bmatrix}\mh_{22}&\mh_{23}\\\mh_{32}&\mh_{33}\end{bmatrix}
=
\begin{bmatrix}\mq_2^\top\mh_{22}&\mq_2^\top\mh_{23}\end{bmatrix}.
\]
Only the columns in $I_3$ are needed: $\mb_3:=\mb(:,I_3)=\mq_2^\top\mh_{23}$. The approximation is $\widetilde{\mh}_{23}=\mq_2\mb_3=\mq_2\mq_2^\top\mh_{23}$, the orthogonal projection of $\mh_{23}$ onto the range of $\mq_2$. The sample and the coefficient factor each cost one block product of width $b$ columns, with $\mh$ and with $\mh^\top$ (which equals $\mh$ in this case). No entries of $\mh_{23}$ are accessed.

The initial factors form a randomized QB approximation~\cite{halkoFindingStructureRandomness2011}. To choose the rank adaptively, we test the current approximation with a fresh Gaussian matrix before adding another block of columns.

For each $i=2,3,\ldots$, draw a new $\mg_{3;i}$ independently of all preceding samples, form $\my_{2;i}=\mh_{23}\mg_{3;i}$ by one block product with $\mh$ as above, and compute the error estimate
\[
\widehat e_i
=\frac{\|\my_{2;i}-\mq_2\mb_3\mg_{3;i}\|_\fro}
{\|\my_{2;i}\|_\fro}
=\frac{\|(\mh_{23}-\widetilde{\mh}_{23})\mg_{3;i}\|_\fro}
{\|\mh_{23}\mg_{3;i}\|_\fro}.
\]
This is the Gaussian estimate \eqref{eq:gaussian_error_indicator} with $\ma=\mh_{23}$, $\widetilde{\ma}=\widetilde{\mh}_{23}$, and $\mOmega=\mg_{3;i}$. \Cref{thm:stability_hutchinson} bounds the true relative error $\|\mh_{23}-\widetilde{\mh}_{23}\|_\fro/\|\mh_{23}\|_\fro$ in terms of this estimate.

If $\widehat e_i\leq\varepsilon$, the test passes and the level is accepted with the current $\widetilde{\mh}_{23}$. Otherwise, the failed sample $\my_{2;i}$ contains directions in the range of $\mh_{23}$ that $\mq_2$ misses, and we reuse it to enlarge the basis. We extend the basis using incremental Householder QR in WY form~\cite{golubMatrixComputations2013}, giving $b$ new orthonormal columns $\mq_{2;i}$ orthogonal to $\mq_2$. Their coefficient rows are $\mb_{3;i}:=\mq_{2;i}^\top\mh_{23}$, obtained by applying $\mh^\top$ to $\mq_{2;i}$ stacked with zeros, exactly as for $i=1$. Append the new columns and rows,
\[
\mq_2\gets\begin{bmatrix}\mq_2&\mq_{2;i}\end{bmatrix},
\qquad
\mb_3\gets\begin{bmatrix}\mb_3\\\mb_{3;i}\end{bmatrix}.
\]
The extension keeps the existing factors and requires one block product with $\mh^\top$ to compute the new coefficients. The next test uses a fresh $\mg_{3;i+1}$, independent of the tested approximation. By symmetry, $\mh_{32}$ has the same relative error with $\widetilde{\mh}_{32}=\widetilde{\mh}_{23}^\top$.

\subsubsection{Level \texorpdfstring{$2$}{2} for a symmetric matrix: inherited errors}
The next level works with the residual $\mr^{(1)}=\mh-\widetilde{\mh}^{(1)}$ of \eqref{eq:peeling_residual}, where $\widetilde{\mh}^{(1)}$ contains the two accepted sibling approximations. The level-$1$ off-diagonal blocks of $\mr^{(1)}$ contain the error $\mh_{23}-\widetilde{\mh}_{23}$ and its transpose, while the diagonal blocks remain unchanged. The level-$2$ targets are the sibling blocks $\mh_{45}$ inside $\mh_{22}$ and $\mh_{67}$ inside $\mh_{33}$. With the shorthand $\mr^{(1)}_{ab}:=\mr^{(1)}(I_a,I_b)$, the residual in level-$2$ blocks is
\[
\mr^{(1)}=
\left[\begin{array}{cc|cc}
\mh_{44}&\mh_{45}&\mr^{(1)}_{46}&\mr^{(1)}_{47}\\
\mh_{54}&\mh_{55}&\mr^{(1)}_{56}&\mr^{(1)}_{57}\\
\hline
\mr^{(1)}_{64}&\mr^{(1)}_{65}&\mh_{66}&\mh_{67}\\
\mr^{(1)}_{74}&\mr^{(1)}_{75}&\mh_{76}&\mh_{77}
\end{array}\right],
\]
where the blocks written with $\mr^{(1)}$ are the four blocks of the level-$1$ error and their transposes. They are controlled by the level-$1$ test, but they are not zero.

We sample both targets using a single block product. Increment $i$ draws fresh standard Gaussian matrices $\mg_{5;i},\mg_{7;i}\in\real^{(n/4)\times b}$, independent of all preceding samples, and places them on the rows in $I_5$ and $I_7$:
\[
\mOmega_i=
\begin{bmatrix}\mathbf 0\\\mg_{5;i}\\\mathbf 0\\\mg_{7;i}\end{bmatrix},
\qquad
\ms_i=\mr^{(1)}\mOmega_i=
\begin{bmatrix}
\mh_{45}\mg_{5;i}+\mr^{(1)}_{47}\mg_{7;i}\\
\mh_{55}\mg_{5;i}+\mr^{(1)}_{57}\mg_{7;i}\\
\mr^{(1)}_{65}\mg_{5;i}+\mh_{67}\mg_{7;i}\\
\mr^{(1)}_{75}\mg_{5;i}+\mh_{77}\mg_{7;i}
\end{bmatrix}.
\]
The rows in $I_4$ sample $\mh_{45}$ and the rows in $I_6$ sample $\mh_{67}$; the other two row blocks are not used. Let $\my_i$ keep the rows of $\ms_i$ in $I_4$ and $I_6$ and set the rows in $I_5$ and $I_7$ to zero. Unlike level $1$, the retained rows are not pure samples of the targets: the terms $\mr^{(1)}_{47}\mg_{7;i}$ and $\mr^{(1)}_{65}\mg_{5;i}$ are errors left by the level-$1$ approximation, and setting the unused rows to zero does not remove them.

For $i=1$, thin QR factorizations of the two nonzero row blocks of $\my_1$ give bases $\mq_4$ and $\mq_6$, each with $b$ orthonormal columns. Stack them in the same $b$ columns, with zeros elsewhere, and compute the coefficients by one block product with $(\mr^{(1)})^\top$:
\[
\mq=\begin{bmatrix}\mq_4\\\mathbf 0\\\mq_6\\\mathbf 0\end{bmatrix},
\qquad
\mb=((\mr^{(1)})^\top\mq)^\top=\mq^\top\mr^{(1)}.
\]
The columns of $\mb$ in $I_5$ and $I_7$ give the two approximations
\[
\begin{aligned}
\mb(:,I_5)&=\mq_4^\top\mh_{45}+\mq_6^\top\mr^{(1)}_{65},
&\qquad
\widetilde{\mh}_{45}&=\mq_4\mb(:,I_5),\\
\mb(:,I_7)&=\mq_4^\top\mr^{(1)}_{47}+\mq_6^\top\mh_{67},
&\qquad
\widetilde{\mh}_{67}&=\mq_6\mb(:,I_7).
\end{aligned}
\]
Thus $\widetilde{\mh}_{45}=\mq_4\mq_4^\top\mh_{45}+\mq_4\mq_6^\top\mr^{(1)}_{65}$. The first term is the projection of $\mh_{45}$ onto the range of $\mq_4$, as at level $1$. The second term is a level-$1$ error that enters through coefficient recovery, because $\mq$ carries both bases in the same columns. The same holds for $\widetilde{\mh}_{67}$, and the transposed blocks follow by symmetry.

For each $i=2,3,\ldots$, draw fresh $\mg_{5;i}$ and $\mg_{7;i}$ and form $\my_i$ by one block product with $\mr^{(1)}$ as above. Let $\widetilde{\mh}^{(2)}$ be the level matrix that holds $\widetilde{\mh}_{45}$, $\widetilde{\mh}_{67}$, and their transposes, with zeros elsewhere. The error estimate of the level is
\[
\widehat e_i=\frac{\|\my_i-\widetilde{\mh}^{(2)}\mOmega_i\|_\fro}{\|\my_i\|_\fro},
\qquad
\widetilde{\mh}^{(2)}\mOmega_i=
\begin{bmatrix}\widetilde{\mh}_{45}\mg_{5;i}\\\mathbf 0\\\widetilde{\mh}_{67}\mg_{7;i}\\\mathbf 0\end{bmatrix}.
\]
This is the Gaussian estimate \eqref{eq:gaussian_error_indicator} for the retained rows and sampled columns of $\mr^{(1)}$:
\[
\ma=\begin{bmatrix}\mh_{45}&\mr^{(1)}_{47}\\\mr^{(1)}_{65}&\mh_{67}\end{bmatrix},
\qquad
\widetilde{\ma}=\begin{bmatrix}\widetilde{\mh}_{45}&\mathbf 0\\\mathbf 0&\widetilde{\mh}_{67}\end{bmatrix},
\qquad
\mOmega=\begin{bmatrix}\mg_{5;i}\\\mg_{7;i}\end{bmatrix}.
\]
The matrix $\ma$ contains the inherited blocks $\mr^{(1)}_{47}$ and $\mr^{(1)}_{65}$, on which $\widetilde{\mh}^{(2)}$ is zero. The stopping test measures more than the compression error of the new sibling blocks.

The test and basis extension follow level $1$, with both bases growing by $b$ columns when a failed sample is reused. As at level $1$, the level is accepted when $\widehat e_i\leq\varepsilon$.

Level $2$ reads only the four blocks of the matrix $\ma$ above: the sibling blocks $\mh_{45}$ and $\mh_{67}$ and the inherited blocks $\mr^{(1)}_{47}$ and $\mr^{(1)}_{65}$. No other block of $\mr^{(1)}$ appears in the samples, in the retained coefficients, or in the test. Level $2$ is therefore adaptive compression of the exact sibling blocks perturbed by these two inherited errors. \Cref{lem:level_reduction} in \S\ref{sec:global_error_analysis} states this for every level and for nonsymmetric matrices.

\subsubsection{General level for a symmetric matrix}
The alternating zero rows in the examples above are described by a row mask. At level $\ell$, the partition splits the rows and the columns into $2^\ell$ blocks. We call these blocks odd or even by their position, counted from $1$, so that the left children are the odd blocks and the right children are the even blocks. Odd rows are the rows in odd blocks, and likewise for even rows and for columns. For $\mx\in\real^{n\times q}$, define
\begin{equation}
\label{eq:peeling_row_mask}
[\mathcal M_\ell(\mx)](I_{2p},:)=\mx(I_{2p},:),
\qquad
[\mathcal M_\ell(\mx)](I_{2p+1},:)=\mathbf 0,
\end{equation}
for each parent $p$ at level $\ell-1$. Thus, $\mathcal M_\ell$ retains odd rows, while $\mathcal M_\ell^{\rm c}(\mx):=\mx-\mathcal M_\ell(\mx)$ retains even rows. Both masks act linearly on rows and are multiplications by diagonal matrices:
$\mathcal M_\ell(\mx)=\mPi_\ell\mx$ and
$\mathcal M_\ell^{\rm c}(\mx)=\mPi_\ell^{\rm c}\mx$, where
$\mPi_\ell:=\mathcal M_\ell(\mi)$ has ones on the odd-block indices and zeros on the even-block indices, and
$\mPi_\ell^{\rm c}:=\mi-\mPi_\ell$.

At level $\ell$, we sample all blocks $\mh_{\alpha\beta}$ with one block product. Increment $i$ draws a fresh standard Gaussian matrix $\mg_i^{(\ell)}\in\real^{n\times b}$, places it on even rows, and keeps the odd rows of the product:
\begin{equation}
\label{eq:peeling_structured_sample}
\mOmega_i^{(\ell)}:=\mathcal M_\ell^{\rm c}(\mg_i^{(\ell)}),
\qquad
\my_i^{(\ell)}:=\mathcal M_\ell\bigl(\mr^{(\ell-1)}\mOmega_i^{(\ell)}\bigr).
\end{equation}
The rows of $\my_i^{(\ell)}$ in $I_\alpha$ equal $\mh_{\alpha\beta}\mg_i^{(\ell)}(I_\beta,:)$ plus the inherited terms from the coarser error blocks in the same rows. Each new Gaussian matrix is independent of the preceding samples. As at level $2$, the mask selects the required rows but does not remove inherited errors from them. In matrix form, $\my_i^{(\ell)}=\mPi_\ell\,\mr^{(\ell-1)}\,\mPi_\ell^{\rm c}\,\mg_i^{(\ell)}$: masking the probe and the product is the same as masking the columns and the rows of the residual. At level $2$, the masked residual $\mPi_2\mr^{(1)}\mPi_2^{\rm c}$ holds the four blocks of the matrix $\ma$ of the level-$2$ example and zeros elsewhere.

For $i=1$, thin QR factorizations of the nonzero row blocks give local bases $\mq_\alpha$, each with $b$ columns. Stack them on odd rows with zeros elsewhere, as at level $2$, to form $\mq$. Since the residual is symmetric, one block product with $\mr^{(\ell-1)}$ gives the coefficients, and the transposed blocks follow by symmetry:
\begin{equation}
\label{eq:peeling_factors_sym}
\mb:=\big(\mr^{(\ell-1)}\mq\big)^\top,
\qquad
\widetilde{\mh}_{\alpha\beta}:=\mq_\alpha\,\mb(:,I_\beta),
\qquad
\widetilde{\mh}_{\beta\alpha}:=\widetilde{\mh}_{\alpha\beta}^\top.
\end{equation}
All sibling factors at the current level share a column count. The level approximation is stored through these factors. As at level $2$, the sample and the coefficients read $\mr^{(\ell-1)}$ on the sibling blocks and on the inherited blocks in the same rows and columns; \Cref{lem:level_reduction} identifies these blocks for every level.

After initialization, each $i=2,3,\ldots$ draws a fresh sample from \eqref{eq:peeling_structured_sample} and computes the error estimate
\begin{equation}
\label{eq:peeling_error_indicator_sym}
\widehat e_i^{(\ell)}
:=
\frac{
\|\my_i^{(\ell)}-\widetilde{\mh}^{(\ell)}\mOmega_i^{(\ell)}\|_\fro
}{
\|\my_i^{(\ell)}\|_\fro
}.
\end{equation}
The test probes the blocks $\mh_{\alpha\beta}$ together with the inherited errors in the same rows; by symmetry, the blocks $\mh_{\beta\alpha}$ have the same relative error, so this is a relative-error test for the full level matrix. 

\subsubsection{Recovery of diagonal blocks}
After $L$ levels have been accepted, the final residual $\mr^{(L)}$ equals $\md_L$ on the leaf diagonal blocks and holds the errors $\mh^{(\ell)}-\widetilde{\mh}^{(\ell)}$ of the accepted levels outside them, by \eqref{eq:hodlr_exact_decomposition}.

By the balanced-tree convention in \S\ref{sec:hodlr}, all leaves have uniform size $s=n/2^L$. We randomize the stacked-identity leaf probe of Martinsson~\cite[\S4]{martinssonCompressingRankStructuredMatrices2016} by assigning independent Rademacher signs to its nonzero entries. This is the block-diagonal specialization of the coloring-based estimator analyzed by Amsel et al.~\cite[\S4.1]{amselFixedSparsityMatrixApproximation2026}, with one color for each column position within a leaf. For each leaf $\tau=2^L,\ldots,2^{L+1}-1$, let $\mPsi_{\tau} := \operatorname{diag} (\psi_{\tau,1},\ldots,\psi_{\tau,s})$, where the $\psi_{\tau,j}$ are independent Rademacher variables ($\pm1$ with probability $1/2$), drawn independently of all samples used in the peeling construction. Let $\boldsymbol\psi\in\real^n$ collect them. The Rademacher variables control the cross-terms in the diagonal recovery error; see \Cref{lem:diagonal_error}. Define
\[
\mP_L
:=
[\,\mPsi_{2^L}\;|\;\mPsi_{2^L+1}\;|\;\cdots\;|\;\mPsi_{2^{L+1}-1}\,]^\top
\in\real^{n\times s}.
\]
\Needspace{4\baselineskip}
With $\ms_L:=\mr^{(L)}\mP_L$ and $\mPsi_{\tau}^2=\mi_s$, we recover the leaf diagonal blocks by
\begin{align}
\widetilde{\md}_L(I_{\tau},I_{\tau})
:=
\ms_L(I_{\tau},:)\mPsi_{\tau} \ =
\md_L(I_{\tau},I_{\tau})
+
\sum_{\tau'\ne\tau}
\mr^{(L)}(I_{\tau},I_{\tau'})
\mPsi_{\tau'}\mPsi_{\tau},
\label{eq:recovered_diagonal_block}
\end{align}
where $\tau'$ runs over the other leaves. Set $\widetilde{\md}_L$ to zero outside the leaf diagonal blocks to complete \eqref{eq:HODLR_approximation}. When $\mh$ is symmetric, the algorithm returns $\tfrac12\bigl(\widetilde{\md}_L+\widetilde{\md}_L^\top\bigr)$ instead of $\widetilde{\md}_L$. Its error is the symmetric part of the original error, so the Frobenius-norm recovery error does not increase, and the returned approximation is symmetric.
\subsubsection{Algorithmic summary and complexity}
\Cref{alg:adaptive_peeling} summarizes the construction. The routine $\texttt{FormSample}(\mr,\ell,b)$ returns $(\mOmega,\my)$ in \eqref{eq:peeling_structured_sample}. The routines $\texttt{Orthogonalize}_\ell$ and $\texttt{ExtendBasis}_\ell$ use incremental Householder QR~\cite{golubMatrixComputations2013} within each nonzero row block, without refactorizing previous samples. They stack bases with zeros elsewhere, and the latter returns only new columns. The routine $\texttt{RecoverDiagonal}$ uses \eqref{eq:recovered_diagonal_block} with Rademacher variables drawn independently of the construction and returns its symmetric part. The rank cutoff $r_\ell^{\rm cut}$ bounds the number of tests at level $\ell$; \Cref{rem:termination} gives its value. If the last test at a level fails, $L$ is not updated, so the level is discarded.
\begin{algorithm}
\caption{Adaptive peeling for a symmetric matrix}
\label{alg:adaptive_peeling}
\begin{algorithmic}[1]
\REQUIRE symmetric $\mh\in\real^{n\times n}$ accessed through matvecs; intended depth $L_{\max}$; relative tolerance $\varepsilon$; block size $b$
\ENSURE HODLR approximation $\widetilde{\mh}$ and accepted depth $L\leq L_{\max}$, with the full-matrix accuracy guarantee of \Cref{thm:global_error} under its stated conditions

\STATE $\mr^{(0)}\gets\mh$;\quad $L\gets0$
\FOR{$\ell=1,\ldots,L_{\max}$}
    \STATE $r_\ell^{\rm cut}\gets$ the rank cutoff of level $\ell$ (\Cref{rem:termination})
    \STATE $[\mOmega_1^{(\ell)},\my_1^{(\ell)}]\gets\texttt{FormSample}(\mr^{(\ell-1)},\ell,b)$
    \STATE $\mq\gets\texttt{Orthogonalize}_\ell(\my_1^{(\ell)})$;\quad $\mb\gets(\mr^{(\ell-1)}\mq)^\top$
    \FOR{$i=2,\ldots,\lceil r_\ell^{\rm cut}/b\rceil+1$}
        \STATE Form the sibling factors of $\widetilde{\mh}^{(\ell)}$ from $\mq$ and $\mb$ as in \eqref{eq:peeling_factors_sym}
        \STATE $[\mOmega_i^{(\ell)},\my_i^{(\ell)}]\gets\texttt{FormSample}(\mr^{(\ell-1)},\ell,b)$
        \STATE $\widehat e_i^{(\ell)}\gets$ the error estimate \eqref{eq:peeling_error_indicator_sym}
        \IF{$\widehat e_i^{(\ell)}\leq\varepsilon$}
            \STATE $\mr^{(\ell)}\gets\mr^{(\ell-1)}-\widetilde{\mh}^{(\ell)}$;\quad $L\gets\ell$;\quad \textbf{break} the inner loop
        \ELSIF{$i\leq\lceil r_\ell^{\rm cut}/b\rceil$}
            \STATE $\mq_i\gets\texttt{ExtendBasis}_\ell(\mq,\my_i^{(\ell)})$;\quad $\mb_i\gets(\mr^{(\ell-1)}\mq_i)^\top$
            \STATE Append the columns $\mq_i$ to $\mq$ and the rows $\mb_i$ to $\mb$
        \ENDIF
    \ENDFOR
    \IF{$\mq$ has at least $r_\ell^{\rm cut}$ columns}
        \STATE \textbf{break} the outer loop \COMMENT{the rank cutoff is reached}
    \ENDIF
\ENDFOR
\STATE $\widetilde{\md}_L\gets\texttt{RecoverDiagonal}(\mr^{(L)},L)$
\RETURN $\widetilde{\mh}=\widetilde{\md}_L+\sum_{\ell=1}^{L}\widetilde{\mh}^{(\ell)}$ and $L$
\end{algorithmic}
\end{algorithm}

We count the cost of \Cref{alg:adaptive_peeling} in three parts: approximation storage, matvecs with $\mh$, and remaining arithmetic. The counts assume a symmetric matrix with no discarded levels. They depend on the accepted ranks $r_1,\ldots,r_L$, where $r_\ell$ is shared by all sibling blocks at level $\ell$, the leaf size $s=n/2^L$, and the block size $b$.

\paragraph{Storage}
The factors of level $\ell$ hold $nr_\ell$ entries. The leaves are symmetric, so one triangle of each leaf is stored, which gives $n(s+1)/2$ entries. The returned approximation therefore stores $n\bigl((s+1)/2+\sum_{\ell=1}^{L}r_\ell\bigr)$ entries. A dense symmetric matrix is also stored by one triangle, about $n^2/2$ entries. \S\ref{sec:results} reports storage in megabytes.

\paragraph{Matvecs}
The construction uses
\begin{equation}
\label{eq:peeling_cost}
N_A=2\sum_{\ell=1}^{L}r_\ell+Lb+s
\end{equation}
matvecs with $\mh$. Each sampled column costs one matvec for the sample and one for its coefficients. Each passing test costs $b$, and diagonal recovery costs one matvec per leaf column. A discarded level adds $2r_{L+1}^{\rm cut}+b=s/2+b$ matvecs if $b$ divides its cutoff.

\paragraph{Arithmetic}
Beyond the matvecs, the algorithm applies the stored levels in the products with $\mr^{(\ell-1)}$ and in the diagonal recovery, orthogonalizes the samples, and applies the current approximation in the tests. Since $b\leq r_\ell$, the total is $\mathcal O\bigl(n\sum_{\ell=1}^{L}r_\ell\sum_{j=1}^{\ell}r_j+ns\sum_{\ell=1}^{L}r_\ell\bigr)$ operations.

\paragraph{Bounded ranks}
Suppose that $r_\ell\leq r$ at every level for a constant $r$ independent of $n$, and that the intended leaf size satisfies $n/2^{L_{\max}}\leq2r$. Then $s\leq2r$, since the algorithm stops at $L=L_{\max}$ or at a level accepted at its cutoff, where $r_L\geq s/2$. Since $L\leq\log_2n$, the approximation requires $\mathcal O(nr\log n)$ storage, and the construction uses $\mathcal O(r\log n)$ matvecs and $\mathcal O(nr^2\log^2n)$ further operations. If a level is discarded, the accepted ranks do not bound $s$, and the bounds apply with the actual leaf size.

The general construction of \Cref{rem:nonsymmetric} stores the factors of both orientations and doubles the sampling, coefficient, and test terms. For the same ranks and leaf size, each count is at most twice the symmetric one. With the cutoff $s_\ell/4$, the same assumptions give $s\leq4r$, and the three orders are unchanged.

\begin{remark}[Nonsymmetric matrices]
\label[remark]{rem:nonsymmetric}
For a nonsymmetric matrix, the blocks $\mh_{\beta\alpha}$ are no longer the transposes of the blocks $\mh_{\alpha\beta}$, and the coefficients require products with $\mh^\top$. We sample the two sibling orientations separately. Label the sample \eqref{eq:peeling_structured_sample} odd--even (OE), after its odd rows and even columns, with $\mOmega_i^{(\ell),\rm OE}:=\mOmega_i^{(\ell)}$ and $\my_i^{(\ell),\rm OE}:=\my_i^{(\ell)}$. The even--odd (EO) sample targets $\mh_{\beta\alpha}$; exchange the masks and use an independent standard Gaussian matrix $\mg_i^{\prime\,(\ell)}$:
\begin{equation}
\label{eq:peeling_structured_sample_rl}
\mOmega_i^{(\ell),\rm EO}:=\mathcal M_\ell(\mg_i^{\prime\,(\ell)}),
\qquad
\my_i^{(\ell),\rm EO}:=\mathcal M_\ell^{\rm c}\bigl(\mr^{(\ell-1)}\mOmega_i^{(\ell),\rm EO}\bigr).
\end{equation}
Thin QR factorizations of the nonzero row blocks of $\my_1^{(\ell),\rm OE}$ and $\my_1^{(\ell),\rm EO}$ give local bases $\mq_\alpha^{\rm OE}$ and $\mq_\beta^{\rm EO}$, stacked on odd and even rows, respectively, into $\mq^{\rm OE}$ and $\mq^{\rm EO}$. One block product with $(\mr^{(\ell-1)})^\top$ per orientation gives
\begin{equation}
\label{eq:peeling_factors}
\begin{aligned}
\mb^{\rm OE}&:=\big((\mr^{(\ell-1)})^\top\mq^{\rm OE}\big)^\top,
&\qquad
\widetilde{\mh}_{\alpha\beta}&:=\mq_\alpha^{\rm OE}\,\mb^{\rm OE}(:,I_\beta),\\
\mb^{\rm EO}&:=\big((\mr^{(\ell-1)})^\top\mq^{\rm EO}\big)^\top,
&\qquad
\widetilde{\mh}_{\beta\alpha}&:=\mq_\beta^{\rm EO}\,\mb^{\rm EO}(:,I_\alpha).
\end{aligned}
\end{equation}
The test combines both orientations in one error estimate,
\begin{equation}
\label{eq:peeling_error_indicator}
\widehat e_i^{(\ell)}
:=
\left(
\frac{
\|\my_i^{(\ell),\rm OE}-\widetilde{\mh}^{(\ell)}\mOmega_i^{(\ell),\rm OE}\|_\fro^2
+
\|\my_i^{(\ell),\rm EO}-\widetilde{\mh}^{(\ell)}\mOmega_i^{(\ell),\rm EO}\|_\fro^2
}{
\|\my_i^{(\ell),\rm OE}\|_\fro^2+\|\my_i^{(\ell),\rm EO}\|_\fro^2
}
\right)^{1/2}.
\end{equation}
With the OE terms alone, \eqref{eq:peeling_error_indicator} reduces to \eqref{eq:peeling_error_indicator_sym}. A failed test extends both bases by $b$ columns, as in the symmetric case, and the sibling factors at a level share a column count. Each increment costs four width-$b$ block products instead of two: one with $\mr^{(\ell-1)}$ and one with $(\mr^{(\ell-1)})^\top$ per orientation. In \Cref{alg:adaptive_peeling}, the sampling, orthogonalization, and coefficient steps then run for both orientations, with the coefficients computed through $(\mr^{(\ell-1)})^\top$ and the test replaced by \eqref{eq:peeling_error_indicator}; the rank cutoff applies to the shared column count.
\end{remark}

\begin{remark}[Termination]
\label[remark]{rem:termination}
Let $s_\ell=n/2^\ell$ be the size of the sibling blocks at level $\ell$. The \emph{rank cutoff} $r_\ell^{\rm cut}$ in \Cref{alg:adaptive_peeling} depends only on $n$ and $\ell$: it is $s_\ell/2$ for the symmetric construction and $s_\ell/4$ for the general construction of \Cref{rem:nonsymmetric}. When the shared column count reaches $r_\ell^{\rm cut}$, level $\ell$ is the last level. One more test accepts or discards it, and the algorithm proceeds to diagonal recovery. Recovery uses the last accepted partition, with larger diagonal blocks if $L<L_{\max}$.

The cutoff balances the matvecs spent on a level against the matvecs it saves in diagonal recovery. In the symmetric construction, compressing a level at rank $r$ costs approximately $2r$ matvecs and reduces the diagonal-recovery width from $2s_\ell$ to $s_\ell$. The balance $2r+s_\ell\approx2s_\ell$ motivates $r_\ell^{\rm cut}=s_\ell/2$. The general construction samples both orientations, so a level at rank $r$ costs approximately $4r$ matvecs, and the same balance gives $r_\ell^{\rm cut}=s_\ell/4$. \end{remark}

\subsection{Global approximation error in adaptive peeling}
\label{sec:global_error_analysis}
In matrix-free peeling, errors in coarser approximations contaminate finer-level samples. This subsection bounds the error of the returned approximation at the accepted depth $L$. The result is \Cref{thm:global_error}: the tolerance bounds the error of the whole matrix, with a constant that does not depend on the depth (matrix size).

We first bound the off-diagonal error when each accepted level satisfies its true local error condition. We then bound the probability of this event and account for diagonal recovery.

\subsubsection{Inherited errors and the level reduction}
For the general analysis, let $\mh\in\real^{n\times n}$ be arbitrary, not necessarily symmetric, and let $\widetilde{\mh}$ denote its HODLR approximation obtained by adaptive peeling with threshold $\varepsilon$. Here, $L$ is the number of accepted levels and may depend on the Gaussian draws. Define
\[
\mk^{(0)}:=0,
\qquad
\mk^{(\ell)}
:=
\mk^{(\ell-1)}
+\mh^{(\ell)}-\widetilde{\mh}^{(\ell)},
\]
so that $\mk^{(\ell)}$ is the sum of the off-diagonal approximation errors of levels 1 through $\ell$, which level $\ell+1$ inherits. It is the off-diagonal part of the residual \eqref{eq:peeling_residual}. To see this, let $\md_\ell$ equal $\mh$ on the diagonal blocks of the level-$\ell$ nodes and zero elsewhere. Applying \eqref{eq:hodlr_exact_decomposition} to the tree truncated at depth $\ell$ gives $\mh=\md_\ell+\sum_{j\leq\ell}\mh^{(j)}$, so $\mr^{(\ell)}=\md_\ell+\mk^{(\ell)}$. The residual equals $\mh$ on these diagonal blocks and $\mk^{(\ell)}$ outside them.

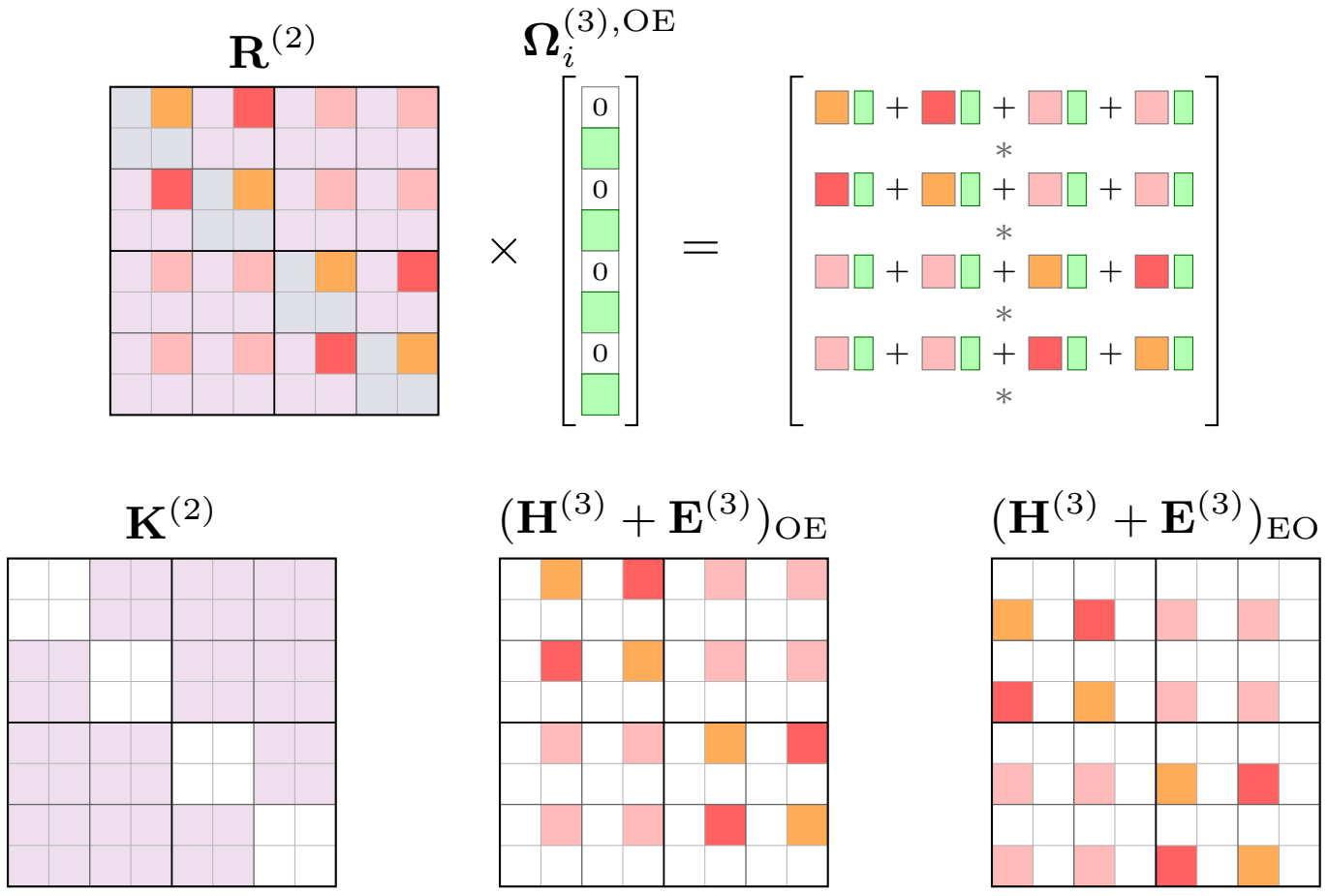
\begin{figure}[ht]
    \centering
    \begin{tikzpicture}[x=.28cm,y=.28cm,
    sib/.style={fill=orange!65},
    eone/.style={fill=red!27},
    etwo/.style={fill=red!62},
    acc/.style={fill=violet!13},
    leaf/.style={fill=blue!13!gray!22},
    gauss/.style={fill=green!30,draw=green!45!black,line width=.2pt},
    zero/.style={fill=white,draw=black!55,line width=.2pt},
    term/.style={draw=black!55,line width=.2pt},
    bracket/.style={draw=black,line width=.4pt}]

\def\fillcell#1#2#3{\path[#1] (#3-1,8-#2) rectangle (#3,9-#2);}

\def\panelgrid{\foreach \k in {1,3,5,7} {\draw[black!30,line width=.15pt] (\k,0)--(\k,8) (0,\k)--(8,\k);}
    \foreach \k in {2,6} {\draw[black!60,line width=.25pt] (\k,0)--(\k,8) (0,\k)--(8,\k);}
    \draw[black,line width=.35pt] (4,0)--(4,8) (0,4)--(8,4);
    \draw[black,line width=.4pt] (0,0) rectangle (8,8);}

\def\errorcell#1#2{\pgfmathtruncatemacro{\sameone}{div(#1-1,4)==div(#2-1,4)}
    \ifnum\sameone=1 \fillcell{etwo}{#1}{#2}\else\fillcell{eone}{#1}{#2}\fi}

\def\panel#1{\foreach \a in {1,...,8} {\foreach \b in {1,...,8} {
        \pgfmathtruncatemacro{\sametwo}{div(\a-1,2)==div(\b-1,2)}
        \pgfmathtruncatemacro{\rsampled}{(mod(\a,2)==0)&&(mod(\b,2)==1)}
        \pgfmathtruncatemacro{\sampled}{(mod(\a,2)==1)&&(mod(\b,2)==0)}
        \ifcase#1\relax
            \ifnum\sametwo=1 \fillcell{leaf}{\a}{\b}\else\fillcell{acc}{\a}{\b}\fi
            \ifnum\sampled=1
                \ifnum\sametwo=1 \fillcell{sib}{\a}{\b}\else\errorcell{\a}{\b}\fi
            \fi
        \or
            \ifnum\sametwo=0 \fillcell{acc}{\a}{\b}\fi
        \or
            \ifnum\sampled=1 \ifnum\sametwo=1 \fillcell{sib}{\a}{\b}\else\errorcell{\a}{\b}\fi\fi
        \or
            \ifnum\rsampled=1 \ifnum\sametwo=1 \fillcell{sib}{\a}{\b}\else\errorcell{\a}{\b}\fi\fi
        \fi
    }}
    \panelgrid}

\begin{scope}[shift={(2.5,0)}]
    \panel{0}
    \node[font=\small] at (4,9) {$\mr^{(2)}$};
    \node at (9.6,4) {$\times$};

\begin{scope}[shift={(11.5,0)}]
        \foreach \a in {1,...,8} {
            \pgfmathtruncatemacro{\isodd}{mod(\a,2)}
            \ifnum\isodd=1
                \path[zero] (0,8-\a) rectangle (.9,9-\a);
                \node[font=\tiny,inner sep=0pt] at (.45,8.5-\a) {$0$};
            \else
                \path[gauss] (0,8-\a) rectangle (.9,9-\a);
            \fi
        }
        \draw[bracket] (-.15,8.25)--(-.45,8.25)--(-.45,-.25)--(-.15,-.25);
        \draw[bracket] (1.05,8.25)--(1.35,8.25)--(1.35,-.25)--(1.05,-.25);
        \node[font=\small] at (.45,9.2) {$\mOmega_i^{(3),\rm OE}$};
        \node at (2.9,4) {$=$};
    \end{scope}

\begin{scope}[shift={(17.2,0)}]
        \foreach \i in {1,...,4} {
            \pgfmathtruncatemacro{\a}{2*\i-1}
            \foreach \j in {1,...,4} {
                \pgfmathtruncatemacro{\b}{2*\j}
                \pgfmathsetmacro{\xl}{2.6*(\j-1)}
                \pgfmathsetmacro{\yb}{8-\a+.1}
                \pgfmathtruncatemacro{\sametwo}{div(\a-1,2)==div(\b-1,2)}
                \pgfmathtruncatemacro{\sameone}{div(\a-1,4)==div(\b-1,4)}
                \ifnum\sametwo=1
                    \path[sib,term] (\xl,\yb) rectangle ++(.8,.8);
                \else\ifnum\sameone=1
                    \path[etwo,term] (\xl,\yb) rectangle ++(.8,.8);
                \else
                    \path[eone,term] (\xl,\yb) rectangle ++(.8,.8);
                \fi\fi
                \path[gauss] (\xl+.95,\yb) rectangle ++(.45,.8);
                \ifnum\j<4 \node[font=\tiny,inner sep=0pt] at (\xl+2,\yb+.4) {$+$};\fi
            }
            \node[font=\scriptsize,inner sep=0pt,text=black!60] at (4.6,8-\a-.55) {$*$};
        }
        \draw[bracket] (-.3,8.25)--(-.6,8.25)--(-.6,-.25)--(-.3,-.25);
        \draw[bracket] (9.5,8.25)--(9.8,8.25)--(9.8,-.25)--(9.5,-.25);
    \end{scope}
\end{scope}

\begin{scope}[shift={(0,-11.5)}]
    \begin{scope}
        \panel{1}
        \node[font=\small] at (4,9) {$\mk^{(2)}$};
    \end{scope}
    \begin{scope}[shift={(12,0)}]
        \panel{2}
        \node[font=\small] at (4,9) {$(\mh^{(3)}+\me^{(3)})_{\rm OE}$};
    \end{scope}
    \begin{scope}[shift={(24,0)}]
        \panel{3}
        \node[font=\small] at (4,9) {$(\mh^{(3)}+\me^{(3)})_{\rm EO}$};
    \end{scope}
\end{scope}
\end{tikzpicture}
     \caption{The level-$3$ reduction in \Cref{lem:level_reduction}. Top: blocks of $\mr^{(2)}$ on the diagonal of the level-$2$ partition contain the original entries of $\mh$. Blue-gray marks the blocks of $\md_2$ outside the OE part, and violet marks the blocks of $\mk^{(2)}$ outside it. The zero probe blocks remove the odd columns of $\mr^{(2)}$ from the product, and the mask discards the even rows (stars). The retained product and the OE coefficients therefore use only the highlighted blocks, those of $\mr^{(2)}_{\rm OE}$: sibling blocks in orange and inherited errors in red. Bottom: the accumulated errors $\mk^{(2)}$ and the OE and EO parts of $\mh^{(3)}+\me^{(3)}$. The OE computation uses the middle panel and the EO computation uses the right panel.}
    \label{fig:error_analysis}
\end{figure}

Because of the masks, only part of $\mk^{(\ell-1)}$ enters the construction at level $\ell$. For an $n\times n$ matrix $\mx$, define its OE part and its EO part at level $\ell$ by
\[
\mx_{\rm OE}:=\mPi_\ell\,\mx\,\mPi_\ell^{\rm c},
\qquad
\mx_{\rm EO}:=\mPi_\ell^{\rm c}\,\mx\,\mPi_\ell,
\]
with the diagonal matrices $\mPi_\ell$ and $\mPi_\ell^{\rm c}$ of \S\ref{sec:adaptivepeeling}. The OE part keeps entries of $\mx$ in odd rows and even columns and is zero elsewhere, and the EO part keeps entries in even rows and odd columns. As a subscript, OE or EO thus selects a part of a matrix. As a superscript, as in \Cref{rem:nonsymmetric}, it labels a quantity computed from the corresponding sample. Since the masks are linear, the two samples \eqref{eq:peeling_structured_sample} and \eqref{eq:peeling_structured_sample_rl} are
\[
\my_i^{(\ell),\rm OE}=\mr^{(\ell-1)}_{\rm OE}\,\mg_i^{(\ell)},
\qquad
\my_i^{(\ell),\rm EO}=\mr^{(\ell-1)}_{\rm EO}\,\mg_i^{\prime\,(\ell)}.
\]
Accordingly, we define the matrix of inherited errors as
\[
\me^{(\ell)}:=\mk^{(\ell-1)}_{\rm OE}+\mk^{(\ell-1)}_{\rm EO}.
\]
The sibling blocks have the same form, $\mh^{(\ell)}=(\md_{\ell-1})_{\rm OE}+(\md_{\ell-1})_{\rm EO}$, because inside a diagonal block of level $\ell-1$ an odd block meets an even block only on the two sibling blocks. Adding the two and using $\mr^{(\ell-1)}=\md_{\ell-1}+\mk^{(\ell-1)}$ gives
\begin{equation}
\label{eq:level_reduction_parts}
\mh^{(\ell)}+\me^{(\ell)}=\mr^{(\ell-1)}_{\rm OE}+\mr^{(\ell-1)}_{\rm EO}.
\end{equation}
For $\ell=2$, written with the node indices of the level-$2$ example in \S\ref{sec:adaptivepeeling}, this matrix is
\[
\mh^{(2)}+\me^{(2)}=
\left[\begin{array}{cc|cc}
\mathbf 0&\mh_{45}&\mathbf 0&\mr^{(1)}_{47}\\
\mh_{54}&\mathbf 0&\mr^{(1)}_{56}&\mathbf 0\\
\hline
\mathbf 0&\mr^{(1)}_{65}&\mathbf 0&\mh_{67}\\
\mr^{(1)}_{74}&\mathbf 0&\mh_{76}&\mathbf 0
\end{array}\right].
\]
The OE part holds the four blocks in the odd rows, those of $I_4$ and $I_6$, and the EO part the other four. For a symmetric matrix, the EO part is the transpose of the OE part. \Cref{fig:error_analysis} shows both parts for $\ell=3$, the first level that inherits errors from two levels. \Cref{lem:level_reduction} shows that level $\ell$ depends on $\mr^{(\ell-1)}$ only through the matrix~\eqref{eq:level_reduction_parts}.

\begin{lemma}[Level-$\ell$ reduction]
\label[lemma]{lem:level_reduction}
At level $\ell$ of adaptive peeling, replace the residual $\mr^{(\ell-1)}$ by $\mh^{(\ell)}+\me^{(\ell)}$. With the same Gaussian draws, the algorithm computes the same level approximation and makes the same stopping decision.
\end{lemma}

\begin{proof}
Write $\mr=\mr^{(\ell-1)}$. The OE sample is $\mr_{\rm OE}\mg_i^{(\ell)}$. The stacked OE basis is supported on odd rows, so $\mq^{\rm OE}=\mPi_\ell\mq^{\rm OE}$, and the factors in \eqref{eq:peeling_factors} use only the even columns of $\mb^{\rm OE}$, with $\mb^{\rm OE}\mPi_\ell^{\rm c}=(\mq^{\rm OE})^\top\mr_{\rm OE}$. The error estimate uses only the samples and factors. Hence the OE computation depends on $\mr$ only through $\mr_{\rm OE}$ and, with the masks exchanged, the EO computation only through $\mr_{\rm EO}$. Since $\mPi_\ell\mPi_\ell^{\rm c}=\mathbf 0$, the OE and EO parts of \eqref{eq:level_reduction_parts} are $\mr_{\rm OE}$ and $\mr_{\rm EO}$. Replacing $\mr$ with this matrix changes neither part, so every increment produces the same factors, test value, and stopping decision. The symmetric construction uses the OE orientation alone.
\end{proof}

Level $\ell$ is therefore an adaptive approximation of the perturbed matrix $\mh^{(\ell)}+\me^{(\ell)}$.

\subsubsection{Local condition and error budgets}
The analysis rests on one condition per level. It compares the level approximation $\widetilde{\mh}^{(\ell)}$ with the perturbed matrix $\mh^{(\ell)}+\me^{(\ell)}$ and not with $\mh^{(\ell)}$. This is the condition that a test can certify, because the samples of level $\ell$ come from the perturbed matrix (\Cref{lem:level_reduction}). It is also sufficient: the inherited errors and the sibling blocks occupy different entries, so the condition leaves a budget for the new compression error (\Cref{lem:local_compression_budget}), and these budgets sum to a bound that does not depend on the depth (\Cref{thm:budgets}).

The sampled test of every level uses the threshold $\varepsilon$. A test that passes certifies a true relative error of at most $3\varepsilon$, by \Cref{thm:stability_hutchinson} with $\eta=1/2$. We therefore define the true error allowance $t:=3\varepsilon$ and assume that $t<1$.

\begin{definition}[Local condition]
\label{def:local_condition}
Level $\ell$ satisfies the \emph{local condition} with allowance $t\in(0,1)$ if
\begin{equation}
\label{eq:adaptive_accuracy}
\|\mh^{(\ell)}+\me^{(\ell)}-\widetilde{\mh}^{(\ell)}\|_\fro
\le
t\,\|\mh^{(\ell)}+\me^{(\ell)}\|_\fro.
\end{equation}
\end{definition}

Let $\mathcal G$ denote the event that every accepted level satisfies the local condition with allowance $t=3\varepsilon$. A level is accepted when its sampled estimate satisfies $\widehat e_i^{(\ell)}\leq\varepsilon$. The estimate is random, so an accepted level can violate \eqref{eq:adaptive_accuracy}, and then $\mathcal G$ does not occur. \Cref{cor:all_levels} bounds the probability of this.

\begin{lemma}[Local compression budget]
\label[lemma]{lem:local_compression_budget}
Suppose that level $\ell$ satisfies \eqref{eq:adaptive_accuracy} with allowance $t\in(0,1)$. Then
\begin{equation}
\label{eq:local_compression_budget}
\begin{aligned}
\|\mh^{(\ell)}-\widetilde{\mh}^{(\ell)}\|_\fro^2
&\leq t^2\|\mh^{(\ell)}+\me^{(\ell)}\|_\fro^2-\|\me^{(\ell)}\|_\fro^2\\
&=t^2\|\mh^{(\ell)}\|_\fro^2-(1-t^2)\|\me^{(\ell)}\|_\fro^2.
\end{aligned}
\end{equation}
\end{lemma}

\begin{proof}
The matrices $\mh^{(\ell)}$ and $\widetilde{\mh}^{(\ell)}$ are supported on the sibling blocks at level $\ell$. The accumulated error $\mk^{(\ell-1)}$ is supported on the union of the coarser-level off-diagonal blocks in \eqref{eq:hodlr_exact_decomposition}. Those blocks are disjoint from the level-$\ell$ sibling blocks. Selecting blocks of $\mk^{(\ell-1)}$ to form $\me^{(\ell)}$ therefore preserves this separation. Consequently,
\begin{align*}
\|\mh^{(\ell)}+\me^{(\ell)}-\widetilde{\mh}^{(\ell)}\|_\fro^2
&=\|\mh^{(\ell)}-\widetilde{\mh}^{(\ell)}\|_\fro^2+\|\me^{(\ell)}\|_\fro^2,\\
\|\mh^{(\ell)}+\me^{(\ell)}\|_\fro^2
&=\|\mh^{(\ell)}\|_\fro^2+\|\me^{(\ell)}\|_\fro^2.
\end{align*}
Applying \eqref{eq:adaptive_accuracy} to these identities gives \eqref{eq:local_compression_budget}.
\end{proof}

Inherited contamination consumes part of the squared budget in \eqref{eq:local_compression_budget}. The identities in the proof show that every approximation supported on the sibling blocks has local error at least $\|\me^{(\ell)}\|_\fro$, with equality for $\widetilde{\mh}^{(\ell)}=\mh^{(\ell)}$. The ratio $\|\me^{(\ell)}\|_\fro/\|\mh^{(\ell)}+\me^{(\ell)}\|_\fro$ is therefore a floor on the relative local error of level $\ell$, the error that remains with exact sibling factors. If the floor equals $t$, the right-hand side of \eqref{eq:local_compression_budget} is zero and only the exact approximation satisfies \eqref{eq:adaptive_accuracy}; if the floor exceeds $t$, no rank satisfies the true local condition. On $\mathcal G$, such a level is not accepted, and the bounds below apply at the smaller depth. A sampled test can still accept it incorrectly; \Cref{thm:certification} bounds the probability of such failures. A floor below $t$ leaves room for compression error. \S\ref{sec:adv-diff-eq} measures the floor on the data-misfit Hessian and shows a level accepted close to its allowance blocking the next one.

Summing the local budgets gives the off-diagonal error at the accepted depth.

\begin{corollary}[Off-diagonal error budgets]
\label{thm:budgets}
On the event $\mathcal G$ that every accepted level satisfies \eqref{eq:adaptive_accuracy}, the off-diagonal approximation error satisfies
\begin{equation}
\label{eq:total_offdiag_error}
\sum_{\ell=1}^{L}\|\mh^{(\ell)}-\widetilde{\mh}^{(\ell)}\|_\fro^2
\leq t^2\sum_{\ell=1}^{L}\|\mh^{(\ell)}\|_\fro^2
=t^2\|\mh-\md_L\|_\fro^2,
\end{equation}
where $\md_L$ is the block-diagonal part of $\mh$ at depth $L$.
\end{corollary}

\begin{proof}
Sum \eqref{eq:local_compression_budget} over the accepted levels. Since $t<1$, each contamination term $-(1-t^2)\|\me^{(\ell)}\|_\fro^2$ is nonpositive. Dropping these terms gives the inequality, and \eqref{eq:hodlr_energy_decomposition} gives the identity.
\end{proof}

The bound \eqref{eq:total_offdiag_error} has no factor that depends on the depth $L$, so the threshold $\varepsilon$ need not be tightened as the number of levels grows.

\subsubsection{Certification of the local conditions}
Every test draws new Gaussian matrices. Once the samples drawn before a test are fixed, \Cref{thm:stability_hutchinson} applies to it and gives a sufficient width $b$ of the test.

\begin{theorem}[Certification by sampled tests]
\label{thm:certification}
Fix the samples drawn before test $i\geq2$ at level $\ell$ of the general construction. Write $\ma=\mh^{(\ell)}+\me^{(\ell)}$ and $\mDelta=\ma-\widetilde{\mh}^{(\ell)}$. If both matrices are nonzero, let $\rho_{\ell,i}$ be the smaller of their stable ranks, and assume
\begin{equation}
\label{eq:peeling_test_width}
b\geq8\rho_{\ell,i}^{-1}\log(2/\delta_{\ell,i}).
\end{equation}
The probability that the test passes and level $\ell$ violates \eqref{eq:adaptive_accuracy} with $t=3\varepsilon$ is at most~$\delta_{\ell,i}$.
\end{theorem}

\begin{proof}
If $\mDelta=0$, the local condition holds directly. If $\ma=0$, \Cref{lem:level_reduction} gives zero coefficients and therefore $\mDelta=0$. Otherwise, $\ma$ and $\mDelta$ are fixed by the earlier samples. By \Cref{lem:level_reduction}, the numerator of the squared error estimate in \eqref{eq:peeling_error_indicator} is the sum of $\|\mDelta\mOmega_i^{(\ell),\rm OE}\|_\fro^2$ and $\|\mDelta\mOmega_i^{(\ell),\rm EO}\|_\fro^2$. The two products have disjoint row supports, so this is $\|\mDelta\mOmega\|_\fro^2$ with $\mOmega=\mOmega_i^{(\ell),\rm OE}+\mOmega_i^{(\ell),\rm EO}$. The denominator is $\|\ma\mOmega\|_\fro^2$ by the same argument. The two probes fill complementary rows with independent standard Gaussian entries. Under \eqref{eq:peeling_test_width}, \Cref{thm:stability_hutchinson} with $\eta=1/2$ gives $\|\mDelta\|_\fro/\|\ma\|_\fro\leq3\widehat e_i^{(\ell)}$ with probability at least $1-\delta_{\ell,i}$. A passing test has $\widehat e_i^{(\ell)}\leq\varepsilon$.
\end{proof}

\begin{corollary}[All accepted levels]
\label{cor:all_levels}
Assign failure allowances $\delta_{\ell,i}$ to the test positions $(\ell,i)$ in advance, including positions that a run may not reach, with $\sum_{\ell,i}\delta_{\ell,i}\leq\delta_{\rm peel}$. If \eqref{eq:peeling_test_width} holds for each reached test with $\ma\ne0$ and $\mDelta\ne0$ as in \Cref{thm:certification},
\begin{equation}
\label{eq:peeling_accuracy_probability}
\mathbb P(\mathcal G^c)\leq\delta_{\rm peel}.
\end{equation}
\end{corollary}

\begin{proof}
Let $\mathcal B_{\ell,i}$ be the event that test $(\ell,i)$ is reached and passes even though level $\ell$ violates \eqref{eq:adaptive_accuracy} with $t=3\varepsilon$. Conditional on the earlier samples that lead to the test, \Cref{thm:certification} bounds this probability by $\delta_{\ell,i}$. Thus $\mathbb P(\mathcal B_{\ell,i})\leq\delta_{\ell,i}$. Since $\mathcal G^c=\bigcup_{\ell,i}\mathcal B_{\ell,i}$, the union bound gives \eqref{eq:peeling_accuracy_probability}.
\end{proof}

For the symmetric construction in \Cref{alg:adaptive_peeling}, the argument applies to the OE part alone. It has the same relative error as the full matrix, the same spectral norms, and half the squared Frobenius norms, so the stable ranks are half as large. Thus the two results hold with $\rho_{\ell,i}/2$ in \eqref{eq:peeling_test_width}.

The union bound allows dependent acceptance decisions. If $T$ bounds the number of possible tests, take $\delta_{\ell,i}=\delta_{\rm peel}/T$ and $b\geq8\log(2T/\delta_{\rm peel})$, because $\rho_{\ell,i}\geq1$ in the general construction and $\rho_{\ell,i}\geq2$ in the symmetric one. The sufficient width grows with $\log T$, while the error constants are independent of depth. This width is $58$ for $\delta_{\rm peel}=0.1$ and $T=64$, a bound on the number of tests when $n=4096$ and $b=32$, so it is about twice the largest block size in the experiments. \S\ref{sec:adv-diff-eq} measures the exact local errors of thirty runs against their thresholds and reports counts at smaller widths.

\subsubsection{Diagonal recovery and the global bound}
Whether $\mathcal G$ occurs is decided by the construction before diagonal recovery, and the Rademacher vector $\boldsymbol\psi$ is drawn independently of the construction samples. In \Cref{lem:diagonal_error} and in the first part of \Cref{thm:global_error}, the construction is therefore fixed, and the expectation and the probability are taken over $\boldsymbol\psi$.

\begin{lemma}[Diagonal-block recovery error]
\label[lemma]{lem:diagonal_error} Fix the off-diagonal construction and consider diagonal-block recovery in \eqref{eq:recovered_diagonal_block} using independent Rademacher variables. Then the expected squared recovery error is the squared off-diagonal error,
\begin{equation}
\label{eq:diagonal_expected_error}
\mathbb E\,\|\md_L-\widetilde{\md}_L\|_\fro^2 = \|\mk^{(L)}\|_\fro^2 = \sum_{\ell=1}^{L}
    \|\mh^{(\ell)}-\widetilde{\mh}^{(\ell)}\|_\fro^2.
\end{equation}
Moreover, for any $0<\gamma<1$, with probability at least $1-\gamma$ over $\boldsymbol\psi$,
\begin{equation}
\label{eq:diagonal_error}
\|\md_L-\widetilde{\md}_L\|_\fro^2 \leq \theta_\gamma \|\mk^{(L)}\|_\fro^2,
\end{equation}
where $\theta_\gamma :=1+\min\Bigl\{\sqrt{2(1-\gamma)/\gamma},\;2\sqrt{\log(1/\gamma)}+2\log(1/\gamma)\Bigr\}$.
\end{lemma}

The squared recovery error is a quadratic form in the Rademacher variables whose mean is the squared off-diagonal error; a variance bound and a subgaussian tail bound give the two branches of $\theta_\gamma$.

\begin{proof}
Let $S=\|\mk^{(L)}\|_\fro^2$ and $X=\|\md_L-\widetilde{\md}_L\|_\fro^2$. By \eqref{eq:recovered_diagonal_block}, the error in each recovered leaf block is the corresponding row block of $\mk^{(L)}\mP_L$ multiplied on the right by $\mPsi_\tau$. Since the leaf diagonal blocks have disjoint supports and $\mPsi_\tau$ is orthogonal, $X=\|\mk^{(L)}\mP_L\|_\fro^2$. Each column of $\mk^{(L)}\mP_L$ is a Rademacher-weighted sum of columns of $\mk^{(L)}$ with the same position within their leaves.

Order the entries of $\boldsymbol\psi$ by their positions within the leaves.
Expanding the squared Euclidean norm of each column and summing
over all columns gives $X=\boldsymbol\psi^\top\mm\boldsymbol\psi$,
where $\mm$ is block diagonal with Gram matrices of the corresponding columns of $\mk^{(L)}$. Hence $\mm$ is positive semidefinite, with $\operatorname{tr}(\mm)=S$, $\|\mm\|_2\leq S$, and $\operatorname{tr}(\mm^2)\leq S^2$. The off-diagonal errors at different levels have disjoint supports, giving $S = \sum_{\ell=1}^{L}
\|\mh^{(\ell)}-\widetilde{\mh}^{(\ell)}\|_\fro^2$. Since $\mathbb E\boldsymbol\psi=\mathbf 0$ and $\mathbb E[\boldsymbol\psi\boldsymbol\psi^\top]=\mi_n$, we get $\mathbb E X=\operatorname{tr}(\mm\mathbb E[\boldsymbol\psi\boldsymbol\psi^\top])=\operatorname{tr}(\mm)=S$.

Using symmetry of $\mm$ and $\psi_j^2=1$, we have $X-\mathbb E X=2\sum_{j<k}\mm_{jk}\psi_j\psi_k$. By independence and zero mean of the Rademacher variables, the products $\psi_j\psi_k$ for distinct unordered pairs are uncorrelated. Hence, $\operatorname{Var}(X)=4\sum_{j<k}\mm_{jk}^2 \leq 2\operatorname{tr}(\mm^2)\leq 2S^2$. For $S>0$, Cantelli's inequality~\cite[Theorem~1.13(a)]{dasguptaProbabilityStatisticsMachine2011} with deviation $S\sqrt{2(1-\gamma)/\gamma}$ gives the first bound in $\theta_\gamma$. For $S=0$, the recovery is exact.

By independence of the entries of $\boldsymbol\psi$ and the inequality $\cosh x\leq e^{x^2/2}$, for any $\boldsymbol\alpha\in\real^n$, we have $\mathbb E\exp(\boldsymbol\alpha^\top\boldsymbol\psi)=\prod_{j=1}^{n}\cosh(\alpha_j)\leq\exp(\|\boldsymbol\alpha\|_2^2/2)$. Since $X=\|\mm^{1/2}\boldsymbol\psi\|_2^2$, the quadratic-form bound of Hsu, Kakade, and Zhang~\cite[Theorem~2.1]{hsuTailInequalityQuadratic2012}, for $u>0$, gives
\[
X\leq S+2\sqrt{\operatorname{tr}(\mm^2)u}+2\|\mm\|_2u
\leq (1+2\sqrt u+2u)S
\]
with probability at least $1-e^{-u}$. Setting $u=\log(1/\gamma)$ gives the second bound in $\theta_\gamma$. Since each
bound holds with probability at least $1-\gamma$, taking the
smaller of the two proves \eqref{eq:diagonal_error}.
\end{proof}

In particular, if the off-diagonal error is zero, the recovery is exact. The second branch of $\theta_\gamma$ grows only logarithmically in $1/\gamma$.

\begin{theorem}[Global approximation error]
\label{thm:global_error}
Let $\widetilde{\mh}$ be the HODLR approximation of $\mh\in\real^{n\times n}$ produced by adaptive peeling at accepted depth $L$ with threshold $\varepsilon<1/3$, followed by diagonal recovery with an independent Rademacher vector $\boldsymbol\psi$, and let $\md_L$ be the block-diagonal part of $\mh$ at depth $L$.
Fix a construction in $\mathcal G$, so every accepted level satisfies \eqref{eq:adaptive_accuracy} with $t=3\varepsilon$. Then the expectation over $\boldsymbol\psi$ satisfies
\begin{equation}
\label{eq:global_expected_error}
\mathbb E\,\|\mh-\widetilde{\mh}\|_\fro
\leq3\sqrt{2}\,\varepsilon\,\|\mh-\md_L\|_\fro,
\end{equation}
and, for any $0<\gamma<1$, with probability at least $1-\gamma$ over $\boldsymbol\psi$,
\begin{equation}
\label{eq:global_error_bound}
\|\mh-\widetilde{\mh}\|_\fro\leq3\,C_\gamma\,\varepsilon\,\|\mh-\md_L\|_\fro,
\qquad
C_\gamma:=(1+\theta_\gamma)^{1/2},
\end{equation}
with $\theta_\gamma$ from \eqref{eq:diagonal_error}. If $\mathbb P(\mathcal G)\geq1-\delta_{\rm peel}$, then \eqref{eq:global_error_bound} holds with probability at least $(1-\delta_{\rm peel})(1-\gamma)$ over the construction samples and $\boldsymbol\psi$.
\end{theorem}

\begin{proof}
Fix a construction in $\mathcal G$. Since the diagonal and off-diagonal errors have disjoint supports,
\[
\|\mh-\widetilde{\mh}\|_\fro^2=\|\mk^{(L)}\|_\fro^2+\|\md_L-\widetilde{\md}_L\|_\fro^2.
\]
By \Cref{lem:diagonal_error}, the expectation of this squared norm over $\boldsymbol\psi$ is $2\|\mk^{(L)}\|_\fro^2$. Jensen's inequality and \Cref{thm:budgets} give
\[
\mathbb E\,\|\mh-\widetilde{\mh}\|_\fro
\leq\Bigl(2\|\mk^{(L)}\|_\fro^2\Bigr)^{1/2}\leq\sqrt{2}\,t\,\|\mh-\md_L\|_\fro
=3\sqrt{2}\,\varepsilon\,\|\mh-\md_L\|_\fro.
\]
With probability at least $1-\gamma$ over $\boldsymbol\psi$, \Cref{lem:diagonal_error} also gives
\[
\|\mh-\widetilde{\mh}\|_\fro^2
\leq(1+\theta_\gamma)\|\mk^{(L)}\|_\fro^2
\leq C_\gamma^2\,t^2\,\|\mh-\md_L\|_\fro^2
=9\,C_\gamma^2\,\varepsilon^2\,\|\mh-\md_L\|_\fro^2,
\]
which proves \eqref{eq:global_error_bound}.

For each construction in $\mathcal G$, the conditional probability of \eqref{eq:global_error_bound} is at least $1-\gamma$. Averaging over the constructions in $\mathcal G$ gives the joint probability bound $(1-\gamma)\mathbb P(\mathcal G)\geq(1-\gamma)(1-\delta_{\rm peel})$.
\end{proof}

For example, $\gamma=10^{-2}$ gives $C_\gamma=3.94$, so that the constant in \eqref{eq:global_error_bound} is $3C_\gamma=11.8$. Since $\|\mh-\md_L\|_\fro\leq\|\mh\|_\fro$, both bounds also hold with $\|\mh\|_\fro$ on the right. The sharper relative bound retains the factor $\|\mh-\md_L\|_\fro/\|\mh\|_\fro$. This factor can be small when the diagonal blocks carry much of the Frobenius norm. The rank also grows in increments of $b$ columns, which can add more columns than the test requires. Both effects help explain why the full-matrix error can be well below the prescribed tolerance. For a symmetric matrix, both bounds also hold for the returned approximation, in which the symmetric part of $\widetilde{\md}_L$ replaces $\widetilde{\md}_L$, because this does not increase the recovery error.

\subsection{Recompression}
\label{sec:recompression}
Adaptive peeling may overestimate the ranks. As noted at the end of \S\ref{sec:global_error_analysis}, the constructed approximation error can be well below the tolerance $\varepsilon$, so smaller ranks suffice. Block size also contributes. The guarantee of \Cref{thm:certification} requires a sufficiently large block size, and large blocks enable efficient BLAS-3 operations. Ranks then grow in large increments, and the last increment can add more columns than its test requires. We therefore strongly recommend recompressing after construction, as is common for hierarchical matrices~\cite[\S6.7]{hackbuschHierarchicalMatricesAlgorithms2015}: we truncate the stored sibling blocks, reducing storage and application cost. There, each block is usually truncated to its own tolerance; here, one error budget is shared by all blocks. Recompression decouples stored ranks from block size: $b$ can be chosen for the tests, with excess columns removed afterward. The recovered diagonal blocks remain unchanged.

First, draw a Gaussian matrix $\mOmega_{\rm val}$ independently of the construction and estimate the relative error by
\[
\widehat e_0:=\frac{\|(\mh-\widetilde{\mh})\mOmega_{\rm val}\|_\fro}
{\|\mh\mOmega_{\rm val}\|_\fro}.
\]
This costs one block product with $\mh$, whose width is the number of columns of $\mOmega_{\rm val}$. We use $20$ columns in all experiments. All remaining truncation steps use the stored factors and are deterministic given these factors and $\widehat e_0$.

Write a stored sibling block of level $\ell$ as $\widetilde{\mh}_{\alpha\beta}=\mq\mb$, where $\mq$ has $r_\ell$ orthonormal columns. If $\mb=\mU\mSigma\mv^\top$ is a thin SVD of $\mb$, then
\[
\widetilde{\mh}_{\alpha\beta}=(\mq\mU)\mSigma\mv^\top
\]
is an SVD of the stored block, because $\mq\mU$ has orthonormal columns. We denote the singular values, the diagonal entries of $\mSigma$, by $\sigma_{\alpha\beta,1}\geq\cdots\geq\sigma_{\alpha\beta,r_\ell}$.

Recompression keeps the leading $k_\ell\leq r_\ell$ singular components of every sibling block of level $\ell$; the ranks $k_\ell$ are chosen below. Let $\widetilde{\mh}_{\rm rec}$ denote the resulting approximation, and $\mDelta:=\widetilde{\mh}-\widetilde{\mh}_{\rm rec}$ denote the discarded part. The blocks occupy disjoint entries, so
\[
\|\mDelta\|_\fro^2
=\sum_{\ell=1}^{L}\sum_{(\alpha,\beta)}\sum_{j>k_\ell}\sigma_{\alpha\beta,j}^2,
\]
where the sum over $(\alpha,\beta)$ includes both directed sibling blocks of every pair at level $\ell$.

The recompressed approximation should still meet the tolerance, $\|\mh-\widetilde{\mh}_{\rm rec}\|_\fro\leq\varepsilon\|\mh\|_\fro$. Since $\mh-\widetilde{\mh}_{\rm rec}=(\mh-\widetilde{\mh})+\mDelta$,
\[
\|\mh-\widetilde{\mh}_{\rm rec}\|_\fro^2
=\|\mh-\widetilde{\mh}\|_\fro^2+\|\mDelta\|_\fro^2
+2\langle\mh-\widetilde{\mh},\mDelta\rangle_\fro,
\]
where $\langle\mx,\my\rangle_\fro:=\tr(\mx^\top\my)$ is the Frobenius inner product. The right-hand side cannot be computed because $\mh$ is available only through matvecs. We approximate it in three steps. First, we omit the last term. The projection part of the construction error is orthogonal to $\mDelta$, but inherited errors in the coefficient factors can contribute to this term. Second, we replace $\|\mh-\widetilde{\mh}\|_\fro$ by its estimate $\widehat e_0\|\mh\|_\fro$. Third, we replace $\|\mh\|_\fro$ by $\|\widetilde{\mh}\|_\fro$, which is computed from the singular values and the dense leaves. By the triangle inequality, the two norms differ by at most $\|\mh-\widetilde{\mh}\|_\fro$, so their relative difference is at most the relative error of the constructed approximation. The condition $\|\mh-\widetilde{\mh}_{\rm rec}\|_\fro^2\leq\varepsilon^2\|\mh\|_\fro^2$ then becomes $\widehat e_0^2\|\widetilde{\mh}\|_\fro^2+\|\mDelta\|_\fro^2\leq\varepsilon^2\|\widetilde{\mh}\|_\fro^2$, that is,
\begin{equation}
\label{eq:rank_pruning_bound}
\|\mDelta\|_\fro^2\leq
(\varepsilon^2-\widehat e_0^2)\|\widetilde{\mh}\|_\fro^2.
\end{equation}
If $\widehat e_0\geq\varepsilon$, nothing is truncated. Because of the three approximations, \eqref{eq:rank_pruning_bound} is a practical budget and does not guarantee a relative error below $\varepsilon$ after recompression.

It remains to choose the ranks $k_1,\ldots,k_L$ within this budget. All sibling blocks at a level keep the same rank, as in construction. Lowering the rank of level $\ell$ from $k_\ell$ to $k_\ell-1$ discards the $k_\ell$th singular component of every sibling block of that level. By the expression for $\|\mDelta\|_\fro^2$ above, this increases $\|\mDelta\|_\fro^2$ by
\[
d_\ell:=\sum_{(\alpha,\beta)}\sigma_{\alpha\beta,k_\ell}^2,
\]
the sum of the squares of the discarded singular values over the level-$\ell$ sibling blocks. The storage saved is the same at every level, $n$ entries for a symmetric matrix, so the best level to lower has the smallest $d_\ell$. Starting from the accepted ranks, $k_\ell=r_\ell$, each step lowers by one the rank of the level with the smallest $d_\ell$. We stop when the next step would violate \eqref{eq:rank_pruning_bound}. Finally, we rewrite each block using the retained singular vectors and values. The work after estimating $\widehat e_0$ is $\mathcal O(n\sum_\ell r_\ell^2)$ for factorization and truncation; each step compares at most $L$ values $d_\ell$.

Although the truncation is deterministic once $\widehat e_0$ is fixed, the ranks still depend on this sampled estimate. We therefore assess the final approximation using a fresh independent Gaussian matrix $\mOmega_{\rm audit}$, also with $20$ columns in our experiments. This audit does not influence rank selection and costs one additional block product with $\mh$. Under the sampling conditions of \Cref{thm:stability_hutchinson}, its error estimate bounds the true relative error with the stated confidence. The experiments report exact errors when the matrix can be formed explicitly.

\section{Numerical results}
\label{sec:results}
We evaluate adaptive peeling on Schur complements arising in PDE solvers (\S\ref{sec:schur_complement}) and on a data-misfit Hessian arising in an inverse problem (\S\ref{sec:adv-diff-eq}). \Cref{tab:notation_results} lists the quantities reported in the tables and figures.

The MATLAB repository with scripts to reproduce all numerical results is available at \url{https://github.com/SMART-Algebra/Adaptive-Matfree-HODLR}.

\begin{table}[ht]
\centering
\caption{Quantities reported in the numerical results.}
\label{tab:notation_results}
\renewcommand{\arraystretch}{1.15}
\begin{tabular}{l p{0.69\textwidth}}
\toprule
$\varepsilon$ & relative tolerance, the threshold of every sampled test \\
$b$ & block size, the number of columns added per increment \\
$L$ & number of accepted levels \\
Recomp.: w/o, with & adaptive peeling without and with the recompression of \S\ref{sec:recompression} \\
Levelwise ranks & shared sibling ranks listed from coarsest to finest \\
Exact error & $\|\mh-\widetilde{\mh}\|_\fro/\|\mh\|_\fro$ against an explicitly formed matrix, affordable for the $4096\times4096$ test matrices \\
$N_A$ & matvecs with the operator, and with its transpose for the general construction, counted during peeling and diagonal recovery \\
Storage & entries stored by the approximation of a symmetric matrix, counted as in \S\ref{sec:adaptivepeeling}, in megabytes at $8$ bytes per entry \\
\bottomrule
\end{tabular}
\end{table}

\subsection{Schur complements arising in PDE solvers}
\label{sec:schur_complement}
The Schur complements of the 2D Poisson problem have constant ranks, so they test the $\mathcal O(r\log n)$ growth of the matvec count derived in \S\ref{sec:adaptivepeeling}, for sizes up to $2^{16}=65{,}536$. The 3D Poisson and convection--diffusion problems, at size $4096$, show the levelwise behavior, including the nonsymmetric case.

\vspace{.5em}
\paragraph{2D problem}
We discretize the 2D Poisson equation on a $51\times n$ interior grid with the five-point stencil, unit spacing, and homogeneous Dirichlet conditions. The middle separator $i_3$ splits the grid into $i_1$ and $i_2$ as shown in \Cref{fig:schur_partition_2D}. Eliminating the two subdomains from the discretization matrix $\ma$ yields the $n\times n$ Schur complement
\begin{equation}
\label{eq:matrixfree_Schur}
\ms =\ma_{33} - \ma_{31} \ma_{11}^{-1} \ma_{13} - \ma_{32} \ma_{22}^{-1} \ma_{23},
\end{equation}
and its matvec uses one solve per subdomain, with sparse Cholesky factors of $\ma_{11}$ and $\ma_{22}$ that are computed once.
\begin{figure}[ht]
\centering
\begin{minipage}[c]{0.42\textwidth}
\centering
\begin{tikzpicture}[x=1.05cm,y=.60cm]
  \foreach \i in {0,...,15} {
    \foreach \j in {0,...,8} {
      \fill[gray!65] (0.2*\i,0.2*\j) circle (0.4pt);
    }
    \fill[blue!70!black] (0.2*\i,0.8) circle (0.8pt);
  }
  \draw[thin] (-0.15,-0.15) rectangle (3.15,1.75);
  \node[fill=white,inner sep=1pt] at (1.5,0.4) {$i_1$};
  \node[fill=white,inner sep=1pt] at (1.5,1.2) {$i_2$};
  \node[blue!70!black,right] at (3.15,0.8) {$i_3$};
  \draw[<->,thin] (-0.15,-0.4) -- (3.15,-0.4) node[midway,below] {$n$};
  \draw[<->,thin] (-0.4,-0.15) -- (-0.4,1.75) node[midway,left] {$51$};
\end{tikzpicture}
 \end{minipage} \begin{minipage}[c]{0.5\textwidth}
\vspace{-1em}
\[
\ma=
\begin{bmatrix}
\ma_{11} & \mathbf{0} & \ma_{13}\\
\mathbf{0} & \ma_{22} & \ma_{23}\\
\ma_{31} & \ma_{32} & \ma_{33}
\end{bmatrix}
\]
\end{minipage}
\caption{Partition of the $51\times n$ grid into the subdomains $i_1$, $i_2$ and the middle separator $i_3$ (left), and the block form of $\ma$ in the order $i_1$, $i_2$, $i_3$ (right).}
\label{fig:schur_partition_2D}
\end{figure}

For these larger matrices we report the audit error $\|(\mh-\widetilde{\mh})\mOmega\|_\fro/\|\mh\mOmega\|_\fro$ using a fresh Gaussian matrix $\mOmega$ with $20$ columns after construction and recompression. \Cref{tab:schur2d_scaling} shows three results. First, adaptive peeling finds bounded ranks: at each tolerance, the accepted ranks are the same for all $n$, consistent with the rank bounds of Chandrasekaran et al.~\cite{chandrasekaran2010numerical}. Recompression reduces them to ranks that just attain the tolerance, with audit errors between $0.76$ and $0.99$ times $\varepsilon$. Second, the cost grows as predicted. With fixed leaf size, each doubling of $n$ adds one level and exactly $2r+b$ matvecs by \eqref{eq:peeling_cost}, so the number of matvecs grows logarithmically with $n$. The storage after recompression grows nearly linearly, by a factor of $2.1$ to $2.2$ per doubling. Third, because the ranks are small, one block of $b=8$ columns already gives an audit error below $5\times10^{-9}$, so the runs at $\varepsilon=10^{-3}$ and $10^{-7}$ are the same before recompression.

\begin{table}[ht]
\centering
\caption{Adaptive peeling for the 2D Poisson Schur complement with block size $b=8$. Without recompression, the rank is the same at every level; with recompression, the range over the levels is given. The storage is after recompression.}
\label{tab:schur2d_scaling}
\setlength{\tabcolsep}{4.6pt}
\renewcommand{\arraystretch}{1.12}
\begin{tabular}{c c c c c c c c c}
\toprule
& & & \multicolumn{2}{c}{\textbf{Levelwise ranks}} & \multicolumn{2}{c}{\textbf{Audit error}} & & \textbf{Storage} \\
\cmidrule(lr){4-5}\cmidrule(lr){6-7}
$\varepsilon$ & $n$ & $L$ & w/o & with & w/o & with & $N_A$ & (MB) \\
\midrule
\multirow{4}{*}{$10^{-3}$} & $2^{13}$ & 9 & 8 & 1--3 & $1.58\times10^{-9}$ & $9.91\times10^{-4}$ & 232 & 1.67 \\
 & $2^{14}$ & 10 & 8 & 1--3 & $2.40\times10^{-9}$ & $8.85\times10^{-4}$ & 256 & 3.60 \\
 & $2^{15}$ & 11 & 8 & 1--3 & $3.30\times10^{-9}$ & $8.75\times10^{-4}$ & 280 & 7.47 \\
 & $2^{16}$ & 12 & 8 & 1--3 & $4.38\times10^{-9}$ & $8.84\times10^{-4}$ & 304 & 15.5 \\
\midrule
\multirow{4}{*}{$10^{-7}$} & $2^{13}$ & 9 & 8 & 5--6 & $1.58\times10^{-9}$ & $8.91\times10^{-8}$ & 232 & 3.77 \\
 & $2^{14}$ & 10 & 8 & 5--6 & $2.40\times10^{-9}$ & $8.89\times10^{-8}$ & 256 & 8.19 \\
 & $2^{15}$ & 11 & 8 & 5--6 & $3.30\times10^{-9}$ & $8.94\times10^{-8}$ & 280 & 17.7 \\
 & $2^{16}$ & 12 & 8 & 5--6 & $4.38\times10^{-9}$ & $9.05\times10^{-8}$ & 304 & 38.0 \\
\midrule
\multirow{4}{*}{$10^{-14}$} & $2^{13}$ & 8 & 16 & 9--10 & $3.51\times10^{-16}$ & $9.27\times10^{-15}$ & 352 & 6.26 \\
 & $2^{14}$ & 9 & 16 & 9--10 & $3.25\times10^{-16}$ & $7.61\times10^{-15}$ & 392 & 13.8 \\
 & $2^{15}$ & 10 & 16 & 9--10 & $3.36\times10^{-16}$ & $8.64\times10^{-15}$ & 432 & 30.0 \\
 & $2^{16}$ & 11 & 16 & 9--10 & $3.40\times10^{-16}$ & $8.90\times10^{-15}$ & 472 & 64.7 \\
\bottomrule
\end{tabular}
 \end{table}

\vspace{.5em}
\paragraph{3D problem}
For the 3D Poisson problem, we use the seven-point stencil on a $51\times64\times64$ interior grid, with unit spacing and homogeneous Dirichlet conditions. The separator is a middle plane, so the Schur complement has size $4096$. We apply \eqref{eq:matrixfree_Schur} using sine-transform subdomain solves.

For the nonsymmetric test, we add $\mathrm{Pe}\,(\partial_xu+\partial_yu+\partial_zu)$ with $\mathrm{Pe}=0.5$ using centered differences on the same grid. Sparse LU factors of the interior blocks give matvecs with the Schur complement and its transpose.

The separator indices are ordered by the kd-tree bisection of \S\ref{sec:hodlr} (\Cref{fig:kdtree_levels}) applied to their $(y,z)$ coordinates, with groups of at most $16$ points, and the Schur complements are approximated in this order.

\begin{table}[ht]
\centering
\caption{Adaptive peeling for the $4096\times4096$ Schur complements of the 3D Poisson problem (symmetric construction) and of the 3D convection--diffusion problem (general construction with both orientations, rank cutoff at half the sibling size), with block size $32$. Storage is in megabytes; for the general construction it counts $2n\sum_{\ell}r_\ell+ns$ entries, for the factors of both orientations and full leaves of size $s$. For comparison, the dense matrix takes $134$~MB, or $67$~MB for one triangle of the symmetric matrix.}
\label{tab:approx_result_schur}
\setlength{\tabcolsep}{2.4pt}
\renewcommand{\arraystretch}{1.12}
\begin{tabular}{l l l l c c c}
\toprule
\textbf{Operator} & $\varepsilon$ & \textbf{Recomp.} & \textbf{Levelwise ranks} & \textbf{Exact error} & \textbf{Storage} & $N_A$ \\
 & & & & & (MB) & \\
\midrule
\multirow{8}{*}{Poisson} & \multirow{2}{*}{$10^{-2}$} & w/o & $[96,\,64,\,64,\,32,\,32,\,32]$ & $8.19\times10^{-4}$ & 11.6 & \multirow{2}{*}{896} \\
 & & \shadecell with & \shadecell $[58,\,32,\,32,\,16,\,16,\,8]$ & \shadecell $9.68\times10^{-3}$ & \shadecell 6.37 & \\
\cmidrule(lr){2-7}
 & \multirow{2}{*}{$10^{-6}$} & w/o & $[224,\,128,\,128,\,64,\,64]$ & $4.84\times10^{-8}$ & 22.0 & \multirow{2}{*}{1504} \\
 & & \shadecell with & \shadecell $[150,\,77,\,79,\,41,\,41]$ & \shadecell $9.76\times10^{-7}$ & \shadecell 14.8 & \\
\cmidrule(lr){2-7}
 & \multirow{2}{*}{$10^{-10}$} & w/o & $[320,\,192,\,192,\,96]$ & $4.99\times10^{-12}$ & 30.4 & \multirow{2}{*}{2144} \\
 & & \shadecell with & \shadecell $[250,\,131,\,131,\,71]$ & \shadecell $9.59\times10^{-11}$ & \shadecell 23.3 & \\
\cmidrule(lr){2-7}
 & \multirow{2}{*}{$10^{-14}$} & w/o & $[448,\,256,\,256]$ & $4.86\times10^{-16}$ & 39.9 & \multirow{2}{*}{2528} \\
 & & \shadecell with & \shadecell $[350,\,184,\,183]$ & \shadecell $9.96\times10^{-15}$ & \shadecell 31.9 & \\
\midrule
\multirow{8}{*}{Conv.--diff.} & \multirow{2}{*}{$10^{-2}$} & w/o & $[96,\,64,\,64,\,32,\,32,\,32]$ & $9.52\times10^{-4}$ & 23.1 & \multirow{2}{*}{1728} \\
 & & \shadecell with & \shadecell $[58,\,32,\,32,\,16,\,16,\,8]$ & \shadecell $9.85\times10^{-3}$ & \shadecell 12.7 & \\
\cmidrule(lr){2-7}
 & \multirow{2}{*}{$10^{-6}$} & w/o & $[224,\,128,\,128,\,64,\,64]$ & $7.07\times10^{-8}$ & 44.0 & \multirow{2}{*}{2880} \\
 & & \shadecell with & \shadecell $[153,\,78,\,80,\,41,\,42]$ & \shadecell $9.71\times10^{-7}$ & \shadecell 30.0 & \\
\cmidrule(lr){2-7}
 & \multirow{2}{*}{$10^{-10}$} & w/o & $[352,\,192,\,192,\,96]$ & $9.96\times10^{-13}$ & 62.9 & \multirow{2}{*}{4160} \\
 & & \shadecell with & \shadecell $[259,\,134,\,136,\,72]$ & \shadecell $9.85\times10^{-11}$ & \shadecell 47.8 & \\
\cmidrule(lr){2-7}
 & \multirow{2}{*}{$10^{-14}$} & w/o & $[448,\,256,\,256]$ & $6.51\times10^{-16}$ & 79.7 & \multirow{2}{*}{4544} \\
 & & \shadecell with & \shadecell $[368,\,192,\,190]$ & \shadecell $9.67\times10^{-15}$ & \shadecell 65.9 & \\
\bottomrule
\end{tabular}
 \end{table}

\Cref{tab:approx_result_schur} shows two results. First, adaptive peeling detects ranks that differ from level to level, for tolerances from $10^{-2}$ to $10^{-14}$ and for both operators. The ranks decrease from the coarsest level to the finest, and the algorithm finds them from matvecs alone, with no rank bounds given in advance. As in 2D, recompression removes the surplus columns of the increments. The exact errors are below the tolerance in all cases, between $0.96$ and $0.996$ times the tolerance after recompression, and recompression reduces the storage by $17\%$ to $45\%$. Second, fewer levels are accepted at tighter tolerances. The required ranks grow as the tolerance tightens, while the sibling size halves from level to level, so the rank cutoff of \Cref{rem:termination}, half the sibling size, is reached at a coarser level: at level $6$ for $\varepsilon=10^{-2}$, at level $5$ for $10^{-6}$ and $10^{-10}$, and at level $3$ for $10^{-14}$. At $\varepsilon=10^{-10}$, the fifth level fails its last test and is discarded.

The convection--diffusion runs use the general construction with the rank cutoff at half the sibling size, as in the symmetric construction, so that both operators reach the same depths. Their accepted ranks equal those of the Poisson runs except at the coarsest level for $\varepsilon=10^{-10}$ ($352$ versus $320$). They use $1.8$ to $1.9$ times as many matvecs, because the two orientations double the sampling, coefficient, and test work while the diagonal recovery is shared.

\Cref{tab:fixed_rank_schur3d} compares adaptive peeling on the Poisson Schur complement with peeling at prescribed ranks and the same depth, without stopping tests. The prescribed ranks are the accepted levelwise ranks $r_\ell$, or one rank for all levels: their maximum, capped at the sibling size, or their rounded mean $\overline r$. With the accepted ranks, the errors are similar, and adaptive peeling uses $96$ to $288$ more matvecs: one width-$32$ stopping test per accepted level and, at $\varepsilon=10^{-10}$, the $160$ matvecs spent on the discarded fifth level. This is the cost of adaptivity. The constant maximum rank meets the tolerances but costs $25\%$ to $45\%$ more matvecs than adaptive peeling, and the constant mean rank misses the tolerances by factors of $4.7$ to $912$. The benefit of adaptivity is that it finds the nonuniform ranks without prior rank information.

\begin{table}[ht]
\centering
\caption{Peeling at prescribed ranks compared with adaptive peeling on the $4096\times4096$ Schur complement of the 3D Poisson problem. All runs use block size $32$ without recompression.}
\label{tab:fixed_rank_schur3d}
\setlength{\tabcolsep}{3pt}
\renewcommand{\arraystretch}{1.12}
\begin{tabular}{c l c c c c}
\toprule
$\varepsilon$ & & \textbf{Adaptive} & \textbf{Fixed, $r_\ell$} & \textbf{Constant, $\max_\ell r_\ell$} & \textbf{Constant, $\overline{r}$} \\
\midrule
\multirow{2}{*}{$10^{-2}$} & $N_A$ & 896 & 704 & 1152 & 700 \\
 & Exact error & $8.19\times10^{-4}$ & $8.10\times10^{-4}$ & $5.10\times10^{-4}$ & $4.70\times10^{-2}$ \\
\midrule
\multirow{2}{*}{$10^{-6}$} & $N_A$ & 1504 & 1344 & 2176 & 1348 \\
 & Exact error & $4.84\times10^{-8}$ & $4.86\times10^{-8}$ & $2.24\times10^{-8}$ & $1.15\times10^{-4}$ \\
\midrule
\multirow{2}{*}{$10^{-10}$} & $N_A$ & 2144 & 1856 & 2688 & 1856 \\
 & Exact error & $4.99\times10^{-12}$ & $4.91\times10^{-12}$ & $4.40\times10^{-12}$ & $9.12\times10^{-8}$ \\
\midrule
\multirow{2}{*}{$10^{-14}$} & $N_A$ & 2528 & 2432 & 3200 & 2432 \\
 & Exact error & $4.86\times10^{-16}$ & $3.30\times10^{-16}$ & $2.58\times10^{-16}$ & $3.55\times10^{-12}$ \\
\bottomrule
\end{tabular}
 \end{table}

The width $32$ of the tests is sufficient for \Cref{thm:certification} at $52$ of the $56$ tests of the Poisson runs in \Cref{tab:approx_result_schur}.{} \Cref{sec:appendix_random_samples} reports ten repetitions of these runs with fresh random numbers.

\subsection{A data-misfit Hessian arising in an inverse problem}
\label{sec:adv-diff-eq}
We consider a linear inverse problem governed by a stationary advection--diffusion--reaction equation. The goal is to infer a spatially varying source term from noisy pointwise observations of the state. We approximate the resulting $4096\times4096$ Hessian. The experiments show how the block size and inherited errors affect the accepted depth, how adaptive peeling compares with a global low-rank approximation, and how the results depend on the random samples.

\begin{figure}[ht]
\centering
\begin{minipage}[c]{0.55\textwidth}
\centering
\includegraphics[width=\linewidth]{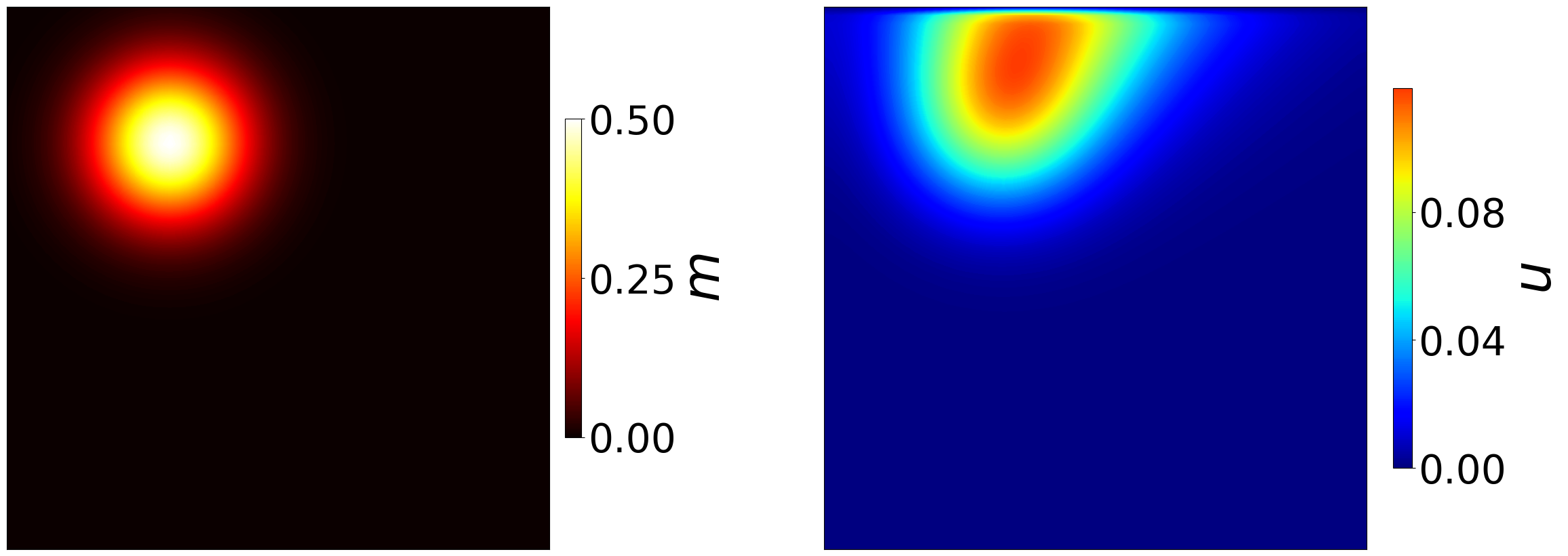}
\end{minipage}\hfill
\begin{minipage}[c]{0.42\textwidth}
\centering
\includegraphics[width=\linewidth]{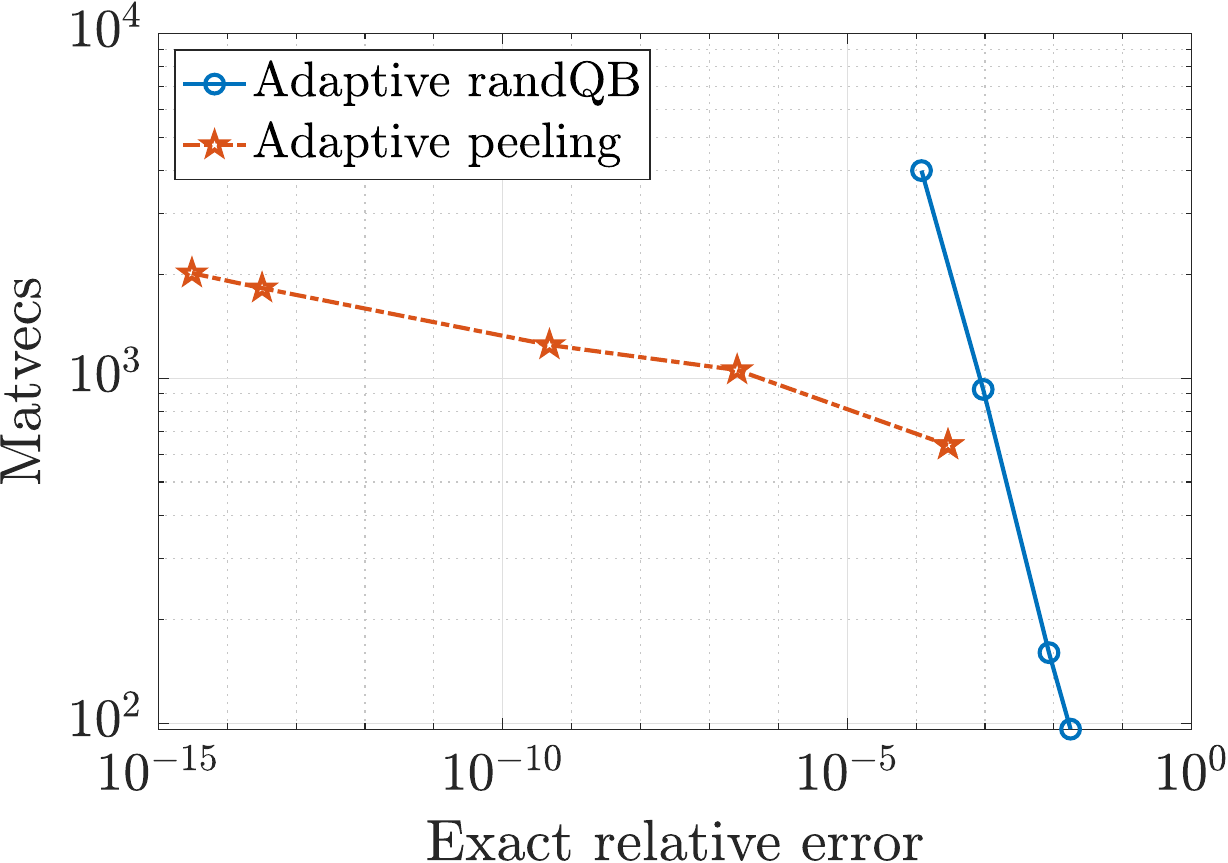}
\end{minipage}
\caption{Left: true source $m_{\mathrm{true}}$ and state $u_{\mathrm{true}}$ for $\kappa=10^{-2}$. Right: adaptive randQB (tolerances $10^{-1}$ to $10^{-4}$) and adaptive peeling (tolerances $10^{-1}$ to $10^{-14}$) for the $4096\times4096$ data-misfit Hessian, both with block size $32$ and without recompression; costs are counted as matvecs, including the stopping tests, and errors are exact.}
\label{fig:hessian_source_state_comparison}
\end{figure}

The state satisfies $-\kappa\Delta u+(x,y)\cdot\nabla u+0.5u=m$ on $(0,1)^2$ with diffusivity $\kappa=10^{-2}$, zero Dirichlet data on the horizontal sides, and zero normal derivative on the vertical sides. We discretize with continuous linear elements on a $64\times64$ nodal grid, splitting each cell along its lower-left to upper-right diagonal. We observe the state at $N_{\mathrm{obs}}=4096$ points on a separate $64\times64$ grid over $[0.001,0.999]^2$. Let $\mf$ be the discrete source-to-observation map. The noise covariance is $\sigma^2\mi$, where $\sigma=10^{-2}\|\mf\vm_{\mathrm{true}}\|_2/\sqrt{N_{\mathrm{obs}}}$, a $1\%$ noise level, and $\vm_{\mathrm{true}}$ is the nodal interpolant of $\tfrac12\exp(-((x-0.3)^2+(y-0.75)^2)/0.02)$. \Cref{fig:hessian_source_state_comparison} (left) shows the source and state.

We approximate the mass-scaled Hessian $\mh_{\mathrm{misfit}}=\sigma^{-2}\ml^{-1}\mf^\top\mf\ml^{-\top}$, where $\ml$ is the lower Cholesky factor of the finite element mass matrix in mesh order. A matvec uses sparse triangular solves and one forward and one adjoint PDE solve. We order the nodes by kd-tree bisection with groups of $16$.

The block size has to be large enough for the fine levels to be compressed. By \Cref{lem:local_compression_budget}, the errors that a level inherits from the coarser levels reduce the budget for its own compression error. \Cref{tab:hessian_blocksize_tol1e-1} shows the effect at $\varepsilon=10^{-1}$. All runs target eight levels, and a level is discarded only if its test still fails when another increment would exceed the sibling-block size. With $b=2$ and $b=4$, the levels are accepted with larger errors, and the first discarded levels, level $7$ for $b=2$ and level $8$ for $b=4$, have true error floors of $0.115$ and $0.119$, above $\varepsilon=0.1$. With $b=8$ and $b=16$, one increment passes at each of the eight levels, and the floors at the finest level are $0.025$ and $0.004$. After recompression, the ranks are $1$ or $2$ at every level for all four block sizes, so the storage is determined by the depth: $0.61$~MB for $b=8$ and $b=16$, compared with $1.36$~MB for $b=2$ and $0.84$~MB for $b=4$. Using $b=16$ reaches the same storage as $b=8$ with $400$ matvecs instead of $208$.

\begin{table}[ht]
\centering
\caption{Hessian block-size comparison at $\varepsilon=10^{-1}$ with the same intended depth of eight for every block size. Levelwise ranks and storage are given without and with recompression; storage is in megabytes, and one triangle of the dense matrix takes $67$~MB.}
\label{tab:hessian_blocksize_tol1e-1}
\setlength{\tabcolsep}{4.2pt}
\renewcommand{\arraystretch}{1.12}
\begin{tabular}{c c r l l rr}
\toprule
 & & & \multicolumn{2}{c}{\textbf{Levelwise ranks}} & \multicolumn{2}{c}{\textbf{Storage (MB)}} \\
\cmidrule(lr){4-5}\cmidrule(lr){6-7}
$b$ & $L$ & $N_A$ & w/o & with & w/o & with \\
\midrule
2 & 6 & 182 & $[4,\,4,\,4,\,2,\,2,\,4]$ & $[1,\,2,\,2,\,2,\,1,\,1]$ & 1.72 & 1.36 \\
4 & 7 & 184 & $[4,\,4,\,4,\,4,\,4,\,4,\,20]$ & $[1,\,2,\,2,\,1,\,1,\,1,\,1]$ & 1.98 & 0.84 \\
8 & 8 & 208 & $[8,\,8,\,8,\,8,\,8,\,8,\,8,\,8]$ & $[1,\,2,\,2,\,1,\,1,\,1,\,1,\,1]$ & 2.38 & 0.61 \\
16 & 8 & 400 & $[16,\,16,\,16,\,16,\,16,\,16,\,16,\,16]$ & $[1,\,2,\,2,\,1,\,1,\,1,\,1,\,1]$ & 4.47 & 0.61 \\
\bottomrule
\end{tabular}
 \end{table}
Global low-rank approximation of data-misfit Hessians is a standard tool in large-scale inverse problems~\cite{flathFastAlgorithmsBayesian2011,bui-thanhComputationalFrameworkInfiniteDimensional2013}. It relies on a fast decay of the Hessian spectrum, which holds when the data inform few parameter modes above the noise. The Hessian here is not of that kind. With $N_{\mathrm{obs}}=n$ observations, $1\%$ noise, and a transport-dominated forward map at $\kappa=10^{-2}$, the data inform most of the parameter space, and the singular values decay slowly. The off-diagonal blocks of the kd-tree ordering are nevertheless compressible, and adaptive peeling meets the tolerances $10^{-6}$ and $10^{-9}$ with five levels of ranks $[128,64,64,64,64]$ and $[160,96,96,64,64]$, respectively.

\Cref{fig:hessian_source_state_comparison} (right) compares adaptive peeling with adaptive randQB~\cite{smithAdaptiveMatrixFreeLowRank2026} applied to the whole Hessian, the method that \S\ref{sec:adaptivepeeling} applies to one sibling block. Adaptive randQB is cheaper at loose tolerances: it reaches an exact error of $1.7\times10^{-2}$ with $96$ matvecs, where peeling spends $640$. Its cost then grows with the rank of the whole Hessian, to $4000$ matvecs at an exact error of $1.2\times10^{-4}$ for a basis of $1984$ columns. Adaptive peeling exploits the hierarchical low-rank structure: it reaches $2.5\times10^{-7}$ with $1056$ matvecs and $3.1\times10^{-15}$ with $2016$.

The Hessian results depend little on the random samples. We repeated the runs with $b=32$ at $\varepsilon=10^{-3}$, $10^{-6}$, and $10^{-9}$ ten times each, with fresh random numbers for the construction, the recompression, the recovery, and the audit; \Cref{tab:sketch_dependence} groups the runs by their accepted ranks. The ranks agree across the ten runs, except in two runs at $\varepsilon=10^{-6}$. There, level $4$ is accepted with an exact local error of $1.19$ and $1.11$ times the threshold, level $5$ is then discarded at the rank cutoff, and the errors after recompression exceed the tolerance by $0.5\%$ and $4.5\%$. In all thirty runs, the error before recompression is below the tolerance.

\begin{table}[ht]
\centering
\caption{Ten runs with fresh random samples for the data-misfit Hessian with block size $32$, grouped by the accepted ranks (number of runs in parentheses). The errors are exact, and a range spans the runs of a group.}
\label{tab:sketch_dependence}
\setlength{\tabcolsep}{6pt}
\begin{tabular}{@{}c l c c@{}}
\toprule
& & \multicolumn{2}{c}{\textbf{Exact error}} \\
\cmidrule(lr){3-4}
$\varepsilon$ & \textbf{Accepted ranks (runs)} & w/o & with \\
\midrule
$10^{-3}$ & $[32,\,32,\,32,\,32,\,32,\,32]$ (10) & $2.56$--$3.13\times10^{-4}$ & $9.17$--$9.31\times10^{-4}$ \\
$10^{-6}$ & $[128,\,64,\,64,\,64,\,64]$ (8) & $2.34$--$3.67\times10^{-7}$ & $9.34$--$9.67\times10^{-7}$ \\
  & $[128,\,64,\,64,\,32]$ (2) & $6.74$--$7.19\times10^{-7}$ & $1.00$--$1.05\times10^{-6}$ \\
$10^{-9}$ & $[160,\,96,\,96,\,64,\,64]$ (10) & $4.50$--$4.88\times10^{-10}$ & $9.02$--$9.92\times10^{-10}$ \\
\bottomrule
\end{tabular}
 \end{table}

The runs also test the analysis. The width $32$ is smaller than the sufficient width of \Cref{thm:certification}, which is between $46$ and $56$ for this matrix. Nevertheless, $156$ of the $158$ accepted levels satisfy \eqref{eq:adaptive_accuracy} with $t=\varepsilon$; the two exceptions are the level-$4$ acceptances above. Narrower tests fail more often: in ten additional runs per block size at $\varepsilon=10^{-1}$, $16$ of $64$ accepted levels exceeded the true threshold with $b=2$, $2$ of $71$ with $b=4$, and none of $80$ with $b=8$.{} \Cref{sec:appendix_random_samples} lists further data for these runs.

\section{Conclusion}
\label{sec:conclusion}
We have developed an adaptive matrix-free randomized method for 
computing an HODLR approximation of a high-dimensional matrix.
Using the proposed method, an HODLR approximation can be computed from matvecs to a prescribed accuracy without knowing the ranks of off-diagonal blocks. Adaptive peeling takes a relative tolerance and finds the ranks level by level with sampled error tests. When the tests use enough random vectors, the relative Frobenius-norm error of the returned approximation is, with high probability, at most a constant times the tolerance, and the constant does not depend on the matrix size. The reason is that each test measures the errors inherited from coarser levels together with the new compression error. Inherited errors therefore cost rank or depth, and not accuracy. In all reported runs, for tolerances from $10^{-1}$ to $10^{-14}$, the errors were below the tolerance before recompression and within $1.1$ times the tolerance after recompression.
 
The analysis bounds the error and not the ranks or the accepted depth. Future work will study how far the accepted ranks are from the smallest ranks that achieve the tolerance, how to choose the tolerances of the levels so that the intended depth is reached at low cost, and the gap between the sufficient and the practical number of random vectors per test.

\section*{Acknowledgments} 
We used Claude (Anthropic) and Codex (OpenAI) to assist with MATLAB code development in the public repository, computational experiments, and the editing and typesetting of this manuscript. Claude suggested the second tail probability bound in \Cref{lem:diagonal_error} and its proof, while the first bound was drafted by the authors. The research ideas, the algorithm, and the design of the analysis and of the experiments are due to the authors. All code, statements, and proofs were checked. The authors assume responsibility for all content.

\appendix
\section{Dependence on the random samples: additional data}
\label{sec:appendix_random_samples}
This appendix extends the experiment on the random samples in \S\ref{sec:adv-diff-eq} to the 3D Poisson Schur complement of \S\ref{sec:schur_complement} and lists further data for the Hessian. We repeated the runs with $b=32$ on both matrices at $\varepsilon=10^{-3}$, $10^{-6}$, and $10^{-9}$ ten times each, with fresh random numbers for the construction, the recompression, the recovery, and the audit. \Cref{tab:sketch_dependence_all} groups the runs by their accepted ranks.

\begin{table}[ht]
\centering
\caption{Ten runs with fresh random samples for the $4096\times4096$ 3D Poisson Schur complement and the data-misfit Hessian at three tolerances, block size $32$, grouped by the accepted ranks (number of runs in parentheses). Recomp.\ columns is the sum of the ranks after recompression over the levels; the errors are exact, and a range spans the runs of the group. The Gaussian matrix for $\widehat e_0$ has $20$ columns.}
\label{tab:sketch_dependence_all}
\setlength{\tabcolsep}{0.9pt}
\renewcommand{\arraystretch}{1.12}
\begin{tabular}{@{}l l c c@{\hspace{6pt}}c@{}}
\toprule
& & \textbf{Recomp.} & \multicolumn{2}{c}{\textbf{Exact error}} \\
\cmidrule(lr){4-5}
\shortstack[l]{\textbf{Problem}\\$\varepsilon$} & \textbf{Accepted ranks (runs)} & \textbf{columns} & w/o & with \\
\midrule
\shortstack[l]{Schur\\$10^{-3}$} & $[128,\,96,\,96,\,64,\,64]$ (6) & 172 & $7.74$--$8.55\times10^{-5}$ & $9.72\times10^{-4}$ \\
  & $[128,\,96,\,96,\,96,\,64]$ (1) & 172 & $8.42\times10^{-5}$ & $9.72\times10^{-4}$ \\
  & $[128,\,96,\,128,\,128]$ (1) & 152 & $7.05\times10^{-5}$ & $9.53\times10^{-4}$ \\
  & $[128,\,96,\,128,\,96,\,64]$ (1) & 172 & $8.54\times10^{-5}$ & $9.72\times10^{-4}$ \\
  & $[128,\,96,\,96]$ (1) & 132 & $5.78\times10^{-5}$ & $9.91\times10^{-4}$ \\
\shortstack[l]{Schur\\$10^{-6}$} & $[224,\,128,\,128,\,64,\,64]$ (10) & 388 & $4.84$--$5.47\times10^{-8}$ & $9.76\times10^{-7}$ \\
\shortstack[l]{Schur\\$10^{-9}$} & $[288,\,160,\,160,\,128]$ (5) & 535 & $6.16$--$6.43\times10^{-11}$ & $9.95\times10^{-10}$ \\
  & $[288,\,160,\,160,\,96]$ (5) & 535 & $6.35$--$6.64\times10^{-11}$ & $9.95\times10^{-10}$ \\
\shortstack[l]{Hessian\\$10^{-3}$} & $[32,\,32,\,32,\,32,\,32,\,32]$ (10) & 39 & $2.56$--$3.13\times10^{-4}$ & $9.17$--$9.31\times10^{-4}$ \\
\shortstack[l]{Hessian\\$10^{-6}$} & $[128,\,64,\,64,\,64,\,64]$ (8) & 241--242 & $2.34$--$3.67\times10^{-7}$ & $9.34$--$9.67\times10^{-7}$ \\
  & $[128,\,64,\,64,\,32]$ (2) & 214--216 & $6.74$--$7.19\times10^{-7}$ & $1.00$--$1.05\times10^{-6}$ \\
\shortstack[l]{Hessian\\$10^{-9}$} & $[160,\,96,\,96,\,64,\,64]$ (10) & 348--349 & $4.50$--$4.88\times10^{-10}$ & $9.02$--$9.92\times10^{-10}$ \\
\bottomrule
\end{tabular}
 \end{table}

The random samples can change the accepted ranks and the accepted depth, and with them the cost. They do not change the accuracy. In all sixty runs, the error before recompression is below the tolerance. After recompression, it exceeds the tolerance only in the two Hessian runs at $\varepsilon=10^{-6}$ described in \S\ref{sec:adv-diff-eq}. These two runs accept level $4$ one increment earlier, and the contamination floors at level $5$ are then $0.95$ and $0.89$ times the threshold.

In three of the six cases, the ten runs accept the same ranks. In the others, some runs differ by one or two increments at one or two levels. For the Schur complement at $\varepsilon=10^{-9}$, the last level is accepted with $96$ columns in five runs and with $128$ in five. At $\varepsilon=10^{-3}$, four runs need more increments at levels $3$ and $4$; one of them accepts level $4$ at the rank cutoff and one discards it, so their depths are four and three instead of five. After recompression, runs of the same depth have the same ranks, up to one column.

The sufficient widths quoted in \S\ref{sec:results} come from the stable-rank quantity $\rho_{\ell,i}/2$ of the symmetric test, computed at every test of one run per tolerance. For the Hessian, it lies between $1.04$ and $1.24$, because one singular value dominates each level matrix; the sufficient width for $\delta_{\rm peel}=0.1$ and $T=64$ is then between $46$ and $56$. For the Poisson runs in \Cref{tab:approx_result_schur} at $\varepsilon=10^{-2}$, $10^{-6}$, and $10^{-10}$, it lies between $1.7$ and $40$, and the width $32$ is sufficient at $52$ of the $56$ tests.

\end{document}